\documentclass[11pt]{article}

\usepackage[latin1]{inputenc}
\usepackage{color}

\usepackage{amsmath}
\usepackage{amssymb}
\usepackage{graphicx}

\usepackage{comment}

\usepackage{accents}

\newcommand{\eps}{\epsilon}
\newcommand{\qed}{$\square$}

\newcommand{\caM}{{\cal M}}

\newcommand{\caF}{{\cal F}}

\newcommand{\al}{\alpha}

\def\lr#1{\langle{#1}\rangle}

\newcommand{\NN}{\mathbb{N}}
\newcommand{\RR}{\mathbb{R}}
\newcommand{\CC}{\mathbb{C}}
\newcommand{\ZZ}{\mathbb{Z}}

\newtheorem{theorem}{Theorem}[section]

\newtheorem{remark}{Remark}[section]

\newtheorem{hypothesis}{Hypothesis}[section]

\newtheorem{definition}{Definition}[section]

\newtheorem{example}{Example}[section]

\newtheorem{proposition}[theorem]{Proposition}

\newtheorem{corollary}[theorem]{Corollary}

\newtheorem{lemma}[theorem]{Lemma}

\numberwithin{equation}{section}

\begin{document}

\begin{center}
\Large
\bf
A Discrete Algorithm for General Weakly Hyperbolic
Systems
\end{center}

\begin{center}
Ferruccio Colombini 
\footnote{Dipartimento di Matematica, Universit\`a di Pisa, Pisa, Italy.
{\tt  ferruccio.colombini@unipi.it}
}
\quad
Tatsuo Nishitani
\footnote{Department of Mathematics,
Graduate School of Science, 
Osaka University, Osaka, Japan. 
{\tt nishitani@math.sci.osaka-u.ac.jp}
}
\quad
 and 
 \quad
 Jeffrey Rauch \footnote{Department of Mathematics, University of Michigan, Ann Arbor, Michigan, USA.
 {\tt rauch@umich.edu}}
\end{center}

\begin{abstract}

This paper studies  the Cauchy problem for  first order systems,
\begin{equation}
\label{eq:cauchy}
Lu \ =\ \partial_tu-\sum_{j=1}^dA_j(t,x)\partial_{x_j}u-B(t,x)u\ =\  f,
\quad
u(0,\cdot) = g\,.
\end{equation}
Assume
 that for $\xi\in \RR^d,$ $\sum A_j\xi_j$ has
only real eigenvalues.    
For coefficients and Cauchy data sufficiently  Gevrey regular
 the Cauchy problem 
has a unique sufficiently  Gevrey regular solution.
We prove  stability and error estimates for the spectral Crank-Nicholson scheme.
Approximate solutions can be computed with accuracy $\eps$ in 
$L^\infty(  [0,T]\times  \RR^d)$
with cost growing at most polynomially in $\eps^{-1}$.  The  proofs uses the 
symmetrizers from \cite{CNR}.
\end{abstract}

\section{Introduction}

\subsection{Hyperbolic background}

Consider  the Cauchy problem 
\eqref{eq:cauchy}.
The coefficients $A_j$ and $B$ are $m\times m$ complex matrix valued
functions that are
 independent of $x$ for $x$ outside a fixed compact set in $\RR^d$. Denote
\[
A(t,x,\xi)
\ :=\
\sum_{j
=
1}^dA_j(t,x)\xi_j\,.
\]
The operator is assumed to satisfy the very weak hyperbolicity condition,
\begin{equation}
\label{eq:Hyperb}
\forall (t,x,\xi)\in \RR\times \RR^d\times \RR^d, \qquad {\rm Spectrum}\, A(t,x,\xi)\subset \RR.
\end{equation}

This hypothesis is best understood by considering first the 
case where the 
coefficients are independent of $t,x$.   In that case,
the initial value problem is solved by Fourier
transform indicated by a hat,
$$
\widehat u(t,\xi) \ =\ 
e^{t(iA(\xi) + B)}\ \widehat g(\xi)\,.
$$
In the case $B=0$, 
the hypothesis implies that for all $\xi \in \RR^d$
$$
\big\|
e^{  iA(\xi) }\ 
\big\|_{ {\rm Hom}(\CC^m)}
\ \lesssim\
\langle   \xi\rangle^{m-1}
\,,
\qquad
\langle\xi\rangle := (1+|\xi|^2)^{1/2}\,.
$$
The Cauchy problem is well set in Sobolev spaces with at worst
a loss of $m-1$ derivatives.
For general $B$ not zero one has the weaker estimate
\begin{equation}
\label{eq:subexp}
\big\|
e^{ iA(\xi) + B }\ 
\big\|_{ {\rm Hom}(\CC^m)}
\ \lesssim\
e^{ c|  \xi|^{(m-1)/m} }
\langle  \xi\rangle^{m-1}
\,.
\end{equation}
This estimate does not allow one to solve the initial value problem for 
$g\in C^\infty_0(\RR^d)$.   However, its subexponential growth shows
that the Cauchy problem is well set in Gevrey spaces.   Those spaces
can be localized by Gevrey partitions of unity so provide a reasonable
setting for the initial value problem
\eqref{eq:cauchy}.
  In the constant coefficient
case, condition 
\eqref{eq:Hyperb} is necessary and sufficient for Gevrey well posedness.

The operators hyperbolic in the sense of 
  Petrowsky and G\aa rding 
\cite{garding}
are characterized by the stronger estimate
\begin{equation}
\label{eq:subexp2}
\big\|
e^{ iA(\xi) + B }\ 
\big\|_{ {\rm Hom}(\CC^m)}
\ \lesssim\
\langle  \xi\rangle^{m-1}
\end{equation}
equivalent to Sobolev solvability with a loss of no more than
$m-1$ derivatives.

Estimate
\eqref{eq:subexp}
  corresponds to 
a sort of  instability at high frequency that  is stronger than permitted  in the constant 
coefficient problems that are hyperbolic in the sense of Petrovsky and G\aa rding
\cite{garding}.

The remarkable fact is that 
provided that the coefficients of $L$ are Gevrey regular,
the Cauchy problem for $L $ is well-posed in  Gevrey classes if and only 
if \eqref{eq:Hyperb} holds.
The sufficiency is a result of  Bronshtein \cite{Bronsh:1}.  
The necessity
is proved in  the trio of articles  \cite{Lax2, mizohata, nishitani}.
 In \eqref{eq:Hyperb},  no hypothesis 
   is made about 
the singularities of  the  characteristic variety of $L$
for $t,x$ fixed, nor on how the geometry of that variety
changes as $t,x$ vary.  The precise Gevrey regularity
required does depend on such structures.  Roughly,
the more variable are the multiplicities the stronger
is the required Gevrey regularity.

\subsection{Algorithm definition}

The present paper
provides more evidence that the weakly
hyperbolic operators  characterized 
by \eqref{eq:Hyperb}
 deserve 
the right to be considered hyperbolic.
We give an algorithm that  computes
approximate solutions
with reasonable computational cost.

Choose $\chi(x)\in C_0^{\infty}(\RR^d)$ with $\chi=1$ in $|x|\leq 2$ and $\chi=0$ for $|x|\geq 2\sqrt{2}$ such that $0\leq \chi\leq 1$. Denote $\chi_h(D)=\chi(h D)$. 
Define  a family of spectral truncations of $G$  by
\begin{equation}
\label{eq:defGh}
G_h(t,x,D)=\chi_{2h} (D) \ \big(iA(t,x,D)+B(t,x)\big)
\chi_{2h}(D),\quad 0<h\leq 1.
\end{equation}
The smoothing operators $G_h$ generate the ordinary differential operators
 $\partial_t-G_h$.
 The resulting  ordinary differential equation is then
approximated by the Crank-Nicholson  scheme.

\begin{definition}\rm 
Define for $n\in \ZZ$,
\begin{equation*}
\label{eq:defGn}
G_h^n(x,D) =
G_h(nk, x, D)
 =
\chi(2hD) \, G(nk, x,D)\, \chi(2hD)\,.
\end{equation*}

The  Crank-Nicholson scheme generating
a sequence $\NN\ni n \mapsto u_h^n$ intended to approximate
$u_h(nk)$ is 
\begin{equation}
\label{eq:CrNi}
\begin{aligned}
\frac{u_h^{n+1}-u_h^{n}}
{k}
\ &=\
G^{n}_h\,
\frac{u_h^{n+1}+u_h^{n}}
{2}
\ +\ 
\chi(2hD)f^{n}\,,
\cr
u^0_h(x) &=\ \chi(2h D)\,g,\qquad f^n\,:=\,f(nk,\cdot).
\end{aligned}
 \end{equation}
\end{definition}

The uniform stability  of the Cauchy problems for $\partial_tu =G_hu$
 is 
proved in \cite{CNR}.  
   This equation has a symmetrizer $R_h=R^*_h$   with $0<c_h<R_h\le 1$.
However as the spectral truncation grows 
the lower bound $c_h$ tends to zero.

Therefore, the straight forward stability arguments that would work for the 
Crank-Nicholson step, as in  \cite{Lax1, YaNo} fail.  The proof of stability
must be at least as hard as the proof of the {\it a priori} estimates
in \cite{CNR}.   Indeed they are more complicated.   The main effort 
follows the strategy in \cite{CNR}.
 We carefully control the additional
errors from discretization in time.   The Crank-Nicholson scheme
is chosen because it is well adapted to estimates using a symmetrizer.

The precise stability result is Theorem \ref{thm:stability}.
The proof  that the approximations converge to the exact
solution is Theorem \ref{thm:convergence}.

For the very special case of operators of the form
$u_{tt} =  a(t) u_{xx}$ with nonnegative Gevery $a$,
the spectral Leap-Frog scheme is analysed in \cite{CR2}.
The computational cost estimates of \cite{CR2} shows that 
the 
 cost of computing with 
error $\eps$ grows no faster than 
polynomially in $\eps^{-1}$.
Virtually identical cost estimates work for our  spectral Crank-Nicholson
scheme.  They are not repeated here.

Constant coefficient problems that 
are hyperbolic in the sense of G\aa rding and Petrowsky
are more strongly  hyperbolic than those studied in this paper.
However variable coefficient operators whose frozen problems
are hyperbolic in this sense need not inherit the Sobolev well
posedness of the constant coefficient problems. 
Stability of difference approximations to constant coefficient
problems hyperbolic in the sense of G\aa rding and Petrowsky
 have
been studied in a number of  works.  We  refer to 
\cite{sabrina} for a review of these.

\section{Main theorems}

\subsection{Definition of the parameter $\theta$}
First we formulate an important property which follows from the assumption \eqref{eq:Hyperb}. Define
\[
{\mathcal H}_r(t,x,\xi,y,\eta;\epsilon)
=
\sum_{|\al+\beta|\leq r}
\frac{\epsilon^{|\al+\beta|}}{\al!\beta!}
D_x^{\al}\partial_{\xi}^{\beta}
A(t,x,\xi)
y^{\al}(-i\eta)^{\beta},\;\;D_{x_j}=-i\partial_{x_j}
\]
then from \cite[Proposition 2.2]{CNR} (see also \cite[(2.3)]{CNR}) it follows that for any compact set $K\subset \RR^d$ and $T>0$ there are $\epsilon_0>0$, $c>0$ such that
\begin{equation}
\label{eq:speccalH}
\zeta\mbox{~is an eigenvalue of~}{\mathcal H}_m(t,x,\xi,y,\eta;\epsilon)
\ \ \Longrightarrow\ \
 |{\mathsf{Im}}\,\zeta|\leq c\,|\epsilon|
\end{equation}
for any $x\in K$, $|\xi|\leq 1$, $|(y,\eta)|\leq 1$, $|\epsilon|\leq \epsilon_0$, $|t|\leq T$.

Following \cite{CNR} introduce an integer $\theta$ defined as follows.

\begin{hypothesis} 
\label{hyp:theta}\rm
The system is  {\it $\theta$-regular} with integer 
$0\le \theta\le m-1$
in the sense that 
 for any $T>0$ and any compact $K\subset \RR^d$ 
there exist  $C>0$, $c>0$ and $\epsilon_0>0$ such that with $N=\max\{2\theta,m\}$
\begin{equation}
\label{eq:Matexp}
\frac{\eps^\theta}
{C\, e^{cs\epsilon}}
 \leq 
 \big\| e^{is{\mathcal H_N(t,x,\xi,y,\eta;\epsilon)}}   \big \|
 \leq
\frac{  C\, e^{cs\epsilon}}
{\epsilon^{\theta}} ,
  \end{equation}
 for all $s\geq 0$, $0<\epsilon\leq \epsilon_0$, $|\xi| \leq 1$, $|(y,\eta)|\leq 1$, $x\in K$, $|t|\leq T$.
 \end{hypothesis}
 \begin{remark}
 \label{rem:deftheta}\rm This definition of $\theta$-regularity is little bit 
 more general than that of \cite[Hypothesis 2.8]{CNR}. 
 Here ${\mathcal H}_r(t,x,\xi,\xi, 0;\epsilon)$ coincides with ${\mathcal H}_r(t,x,\xi; \epsilon)$ in \cite{CNR}.
 \end{remark}
 \begin{example}
 \label{examp:9and10}\rm 
 When \eqref{eq:Hyperb} holds,
  Hypothesis  \ref{hyp:theta} always holds with $\theta=m-1$. If $A(t,x,\xi)$ is uniformly diagonalizable then Hypothesis  \ref{hyp:theta} holds with $\theta=0$ (for the proof see \cite[Examples 2.9 and 2.10]{CNR}).
 \end{example}
 \begin{example}
 \label{examp:11}\rm Suppose \eqref{eq:Hyperb}. Assume that there exists $T=T(t,x,\xi,y,\eta;\epsilon)$ with bounds on $\|T\|$ and $\|T^{-1}\|$ independent of $(t,x,\xi,y,\eta;\epsilon)$ such that $T^{-1}{\mathcal H}_mT$ is a direct sum $\oplus A_j$ where the size of $A_j$ is at most $\mu$. Then Hypothesis  \ref{hyp:theta} holds with $\theta=\mu-1$ (for the proof see \cite[Example 2.11]{CNR}).
 \end{example}
 %
 
 \subsection{Recall the continuous case}

Let
\[
G(t,x,D)=iA(t,x,D)+B(t,x)
\]
then $Lu=f$ is written
\[
\partial_t u=Gu+f.
\]
Denote
\begin{equation}
\label{eq:defanglexi}
\lr{\xi}_{\ell}=
\sqrt{\ell^2+|\xi|^2}
\ =\ 
\ell\sqrt{1\, +\, |\xi/\ell|^2}
\end{equation}
where $\ell\geq 1$ is a positive  parameter. When $\ell=1$ we omit the suffix $\ell$ and write $\lr{\xi}_1=\lr{\xi}$. 
\begin{definition}
\label{dfn:gs}
If $1<s<\infty$, the function $a(x)\in  C^{\infty}(\RR^d)$ belongs  to  $G^{s}(\RR^d)$ if  there exist $C>0$, $A>0$ such that
\[
\forall  x\in {\mathbb R}^d,\quad
\forall  \alpha\in\NN^d,
\qquad
|\partial_x^{\al}a(x)|
\ \leq\ C A^{|\al|}|\al|!^s\,.
\]
\end{definition}

Recall \cite[Proposition 4.4]{CNR}.
 \begin{proposition}
 \label{pro:CNR:1}Suppose Hypothesis \ref{hyp:theta} is satisfied.
Define
$$
s=
    \frac{ 2+6\theta }
    {1+6\theta},\qquad \rho=\frac{1}{s}\,,
    \qquad
    \nu:= \theta(1-\rho)\,.
    $$
For some $1<s'\le s$ suppose that 
$A_j(t,x)$ {\rm(}resp. $B(t,x)${\rm)} are lipschitzian
{\rm(}resp. continuous{\rm)} in time uniformly on compact sets with 
values in $G^{s'}(\RR^d)$. Then there exist $T>0$, ${\hat c}>0$, $C>0$ and $\ell_0>0$ such that  for all $u$ such  that $e^{(T-{\hat c}\,t)\lr{D}_{\ell}^{\rho}}\partial_{t,x}^{\gamma}u\in L^1([0,T/{\hat c}];H^{\nu}(\RR^d))$ for $|\gamma|\leq 1$ one has
 \begin{equation}
 \label{eq:conti_estimate}
 \begin{split}
 \|\lr{D}^{-\nu}_{\ell}e^{(T-{\hat c}\,t)\lr{D}_{\ell}^{\rho}}u\|^2\leq C\|\lr{D}_{\ell}^{\nu}e^{T\lr{D}_{\ell}^{\rho}}u(0)\|^2\\
 +C\int_0^t\|\lr{D}_{\ell}^{-\nu}e^{(T-{\hat c}\,t')\lr{D}_{\ell}^{\rho}}(\partial_t-G)u(t')\|^2dt'
 \end{split}
 \end{equation}
 for $0\leq t\leq T/{\hat c}$ and $\ell\geq \ell_0$.
 \end{proposition}
  This is a small improvement  of 
 \cite[Proposition 4.4]{CNR}.  
 Here is a 
 sketch of the easy proof: As noted in Remark \ref{rem:deftheta} we use ${\mathcal H}_r(t,x,\xi,y,\eta; \epsilon)$ instead of ${\mathcal H}_r(t,x,\xi;\epsilon)$ in \cite{CNR} and make the same choice \eqref{eq:choice} below for $s, \epsilon, \xi,y,\eta$ where $\chi_h\equiv 1$, $\chi_{2h}\equiv 1$ and ${\bar \tau}-\tau=T-at$. Therefore \eqref{eq:expH} below holds for $0\leq T-at\leq {\bar \tau}$ which gets rid off the constraint $T-at\geq c$ with some $c>0$ that we have assumed in \cite{CNR}. This enables us to take $T_1=T$ in \cite[Proposition 4.4]{CNR}. In the estimate \eqref{eq:conti_estimate} the weight for $Lu$ is   improved from $\lr{D}_{\ell}^{3\nu}$ to $\lr{D}_{\ell}^{-\nu}$.  That   proof is also easy.

 \begin{corollary}
 \label{cor:fixell}There exist $T>0$,  ${\hat c}>0$, $C>0$ and $\ell_0>0$  such that for all $u$ satisfying $\partial_tu=Gu$ one has
 \begin{equation}
 \label{eq:conti_estimate:2}
 \|\lr{D}^{-\nu}_{\ell}e^{(T-{\hat c}\,t)\lr{D}^{\rho}_{\ell}}u\|\leq C\|\lr{D}^{\nu}_{\ell}e^{T\lr{D}^{\rho}_{\ell}}u(0)\|
 \end{equation}
 for $0\leq t\leq T/{\hat c}$  and $\ell\geq \ell_0$.
 \end{corollary}
The proof of \cite[Theorem 1.3]{CNR} gives
 \begin{proposition}
 \label{pro:conti_existence}Assume the same assumption as in Proposition \ref{pro:CNR:1} and $e^{T\lr{D}^{\rho}}g\in H^{\nu}(\RR^d)$. Then there exists
 a unique  $u$ satisfying
 \[
 \partial_tu=Gu\,, \qquad t\in (0,T/{\hat c}),\quad u(0,\cdot)=g
 \]
 such that $e^{(T-{\hat c}\,t)\lr{D}^{\rho}}u\in L^{\infty}([0,T/{\hat c}];H^{-\nu}(\RR^d))$.
 \end{proposition}
 %

\subsection{Stability and error estimates}

The Crank-Nicholson scheme defined in 
  \eqref{eq:CrNi} is equivalent to 
 \begin{equation}
 \label{eq:schemeout}
 \big(
 I - \frac{k}{2}\, G^{n}_h\big)
 u_h^{n+1}
 \ =\ 
  \big(
 I +\frac{k}{2}\, G^{n}_h\big)u_h^n
\ +\ 
k\,\chi_{2h}\,f^{n}
 \,.
 \end{equation}
 Note that
 \begin{align*}
\big\|\frac{k}{2}\,G^n_h u \big\|\leq \frac{k}{2}\,\big\| \lr{D}\chi_{2h}\lr{D}^{-1}G(nk,x,D)\chi_{2h}u\big\|\\
 \leq \frac{\sqrt{3}}{2}k\,h^{-1}\big\| \lr{D}^{-1}G(nk,x,D)\chi_{2h}u\big\|  \leq {\bar C} k\,h^{-1}\|u\|
 \end{align*}
 where
 \[
 {\bar C}=\frac{\sqrt{3}}{2}\sup_{0\leq t\leq T}\big\| \lr{D}^{-1}G(t,x,D)\big\|_{{\mathcal L}(L^2,L^2)}.
 \]
Assuming ${\bar C}\,k\,h^{-1}<1$ one has
    \begin{equation}
  \label{eq:Neumann}
 \big(
 I - \frac{k}{2}\, G^{n}_h\big)^{-1}
 \ =\ 
 \sum_{j=0}^\infty
 \Big(
 \frac{k}{2}\, G^{n}_h
 \Big)^j\,,
 \end{equation}
 and
  $u^{n+1}_h$ is given by 
   $$
 u^{n+1}_h
 =
 \big(
 I - \frac{k}{2}\, G^{n}_h\big)^{-1}
  \Big(
  \big(
 I +\frac{k}{2}\, G^{n}_h\big)u^n_h
\ +\ 
k\,\chi_{2h}\,f^{n}
 \Big)\,.
 $$

Reasoning term by term in \eqref{eq:Neumann}, 
$\big( I - \frac{k}{2}\, G^{n}_h\big)^{-1}$
 maps functions with spectrum in ${\rm supp}\,\chi_{2h}(\cdot)$
 to themselves.  Therefore,
 \begin{equation}
\label{eq:smallspec}
{\rm supp}\, \caF({u^n_h})
\ \subset\
 {\rm supp}\,\chi_{2h}(\cdot)
\,.
\end{equation}
 \begin{theorem}
 \label{thm:stability}  
 Make the same assumption as in Proposition \ref{pro:CNR:1}. Then 
there exist ${\bar \tau}>0$, ${\bar \beta}>0$, ${\bar a}>0$, ${\bar h}>0$ and $C>0$ such that the estimate
 \begin{align*}
\|\lr{D}^{-\nu} e^{({\bar\tau}-{\bar a}t_n)\lr{D}^{\rho}} u^n_h\|^2\leq C\Big(\|\lr{D}^{\nu}e^{{\bar\tau}\lr{D}^{\rho}}g\|^2
  +k\,\sum_{j=0}^{n-1}\|\lr{D}^{-\nu}e^{({\bar\tau}-{\bar a}t_j)\lr{D}^{\rho}}f^j\|^2\Big)\\
  \leq C\Big(\|\lr{D}^{\nu}e^{{\bar\tau}\lr{D}^{\rho}}g\|^2
  +\sup_{0\leq j\leq n-1}\|\lr{D}^{-\nu}e^{({\bar\tau}-{\bar a}t_j)\lr{D}^{\rho}}f^j\|^2\Big)
  \end{align*}
holds   for any $n\in \NN$, $k>0, h>0$ satisfying $t_n=nk\leq {\bar \tau}/{\bar a}$, $kh^{-1}\leq {\bar\beta}$ and $0<h\leq {\bar h}$ where $\nu=\theta(1-\rho)$.
  \end{theorem}
 A more precise estimate of the stability is given in Proposition \ref{pro:CrNi_stability}.
\begin{theorem}
\label{thm:convergence}
In addition to  the  assumption  in Proposition \ref{pro:CNR:1}, assume that $A_j(t,x)$  and $B(t,x)$ are $C^1$
 in time uniformly on compact sets with 
values in $G^{s'}(\RR^d)$. Then 
there exist ${\bar \tau}>0$, ${\bar \beta}>0$, ${\bar a}>0$, ${\bar h}>0$ and $C>0$ such that for 
an  exact solution $u$ to \eqref{eq:cauchy} with Cauchy data  $g$ satisfying  $\lr{D}^{2+\nu}e^{{\bar\tau}\lr{D}^{\rho}}g\in L^2$ one has
\begin{align*}
\|\lr{D}^{-\nu}e^{({\bar\tau}-{\bar a}t_n)\lr{D}^{\rho}}(u(t_n)-u^n_h)\|
\leq C\,(k+h)\|\lr{D}^{2+\nu}e^{{\bar\tau}\lr{D}^{\rho}}g\|
\end{align*}
and
\[
\|e^{({\bar\tau}-{\bar a}t_n)\lr{D}^{\rho}}(u(t_n)-u^n_h)\|
\leq C\,(k+h)h^{-\nu}\|\lr{D}^{2+\nu}e^{{\bar\tau}\lr{D}^{\rho}}g\|
\]
for any $n\in \NN$, $k>0, h>0$ satisfying $t_n=nk\leq {\bar \tau}/{\bar a}$, $kh^{-1}\leq {\bar\beta}$ and $0<h\leq {\bar h}$.
\end{theorem}
\begin{corollary}
\label{cor:CrNi_convergence} With  the same assumptions as in Theorem \ref{thm:convergence} 
there exist ${\bar \tau}>0$, ${\bar \beta}>0$, ${\bar a}>0$, ${\bar h}>0$ and $C>0$ such that for an  exact solution $u$ to \eqref{eq:cauchy} with Cauchy data  $g$ satisfying $\lr{D}^{2+\nu}e^{{\bar\tau}\lr{D}^{\rho}}g\in L^2$ one has
\[
\|u(t_n)-u^n_h\|
\leq C\,(k+h)h^{-\nu}\|\lr{D}^{2+\nu}e^{{\bar\tau}\lr{D}^{\rho}}g\|
\]
for any $n\in \NN$, $k>0, h>0$ satisfying $t_n=nk\leq {\bar \tau}/{\bar a}$, $kh^{-1}\leq {\bar\beta}$ and $0<h\leq {\bar h}$.
\end{corollary}
 \begin{remark}
 \label{rem:rhoandnu}
 \rm Note that
 \begin{equation}
 \label{eq:rhoandnu}
 \rho\geq \frac{1+6\theta}{2+6\theta}
 \quad
 \Longleftrightarrow
 \quad \rho\geq 3\nu +\frac{1}{2}
 \end{equation}
so that one has $\rho\geq 3\nu+1/2$ under the assumption of Theorems \ref{thm:stability} and \ref{thm:convergence}.
 \end{remark}

\section{Stability  for the spectral Crank-Nicholson scheme}

\subsection{Spectral truncated weight for Crank-Nicholson scheme}

Taking \eqref{eq:CrNi} into account define
 spectral truncated weights $W_h(t,D)$  by
\[
W_h(t,\xi)
\ :=\
e^{(T-t)\lr{\xi}_{\ell}^{\rho}\chi_h(\xi)}
\]
and for $n\in\NN$
\[
W_h^n(\xi)\ :=\
W_h(ank,\xi)
\]
where $a>0$ is a positive parameter which will be fixed later. In what follows we always assume that the parameters $a>0, k>0, \ell>0, h>0$ are constrained to satisfy
\begin{align}
\label{eq:constraint:1}
0<h\leq \ell^{-1},\quad k\,h^{-1}\leq 1/2\,{\bar C},\quad a\,k\,h^{-\rho} \leq \log{2}/3.
\end{align}
Since
 $a\lr{\xi}_{\ell}^{\rho}\chi_h\leq 3\,a\, h^{-\rho}$
because $\lr{\xi}_{\ell}\leq 3h^{-1}$ if $\chi_h(\xi)\neq 0$,
it follows that 
\begin{equation}
\label{eq:W_h}
1/2\leq e^{-ak\lr{\xi}_{\ell}^{\rho}\chi_h}\leq 1.
\end{equation}
Here recall \cite[Definition 2.3]{CNR}.

\begin{definition}
\label{dfn:rhodelta} \rm
For $0<\delta\le \rho\le 1$, the family   $a(x,\xi;\ell)\in C^{\infty}(\RR^d\times\RR^d)$ 
indexed by $\ell$
belongs to  ${\tilde S}^m_{\rho,\delta}$ if for all $\al$, $\beta\in \NN^d$ there is $C_{\al\beta}$ independent of 
$\ell\ge 1,x,\xi$ such that
\[
\big|
\partial_x^{\beta}\partial_{\xi}^{\al}a(x,\xi;\ell)
\big|
\ \leq\
 C_{\al\beta}\
 \lr{\xi}_\ell^{m-\rho|\al|+\delta|\beta|}.
\]
Denote  ${\tilde S}^m={\tilde S}^m_{1,0}$.
\end{definition}
Since $|\partial_{\xi}^{\alpha}\chi_h|\leq C_{\alpha}h^{|\alpha|}$ and $2h^{-1}\leq \lr{\xi}_{\ell}\leq 3h^{-1}$ on the support of $\partial_{\xi}^{\alpha}\chi_h$ for $|\alpha|\geq 1$ it is clear that $\chi_h\in {\tilde S}^0$.

We examine to  what extent $W^n_h$ satisfies the Crank-Nicholson scheme \eqref{eq:CrNi}. 
\begin{lemma}
\label{lem:hyoome}Assume \eqref{eq:constraint:1} then one can write 
\begin{equation}
\label{eq:W_CrNi}
\frac{W_h^{n+1}-W^n_h}{k}=-2\,a\,\omega_h\,\chi_h \frac{W^{n+1}_h+W^n_h}{2}
\end{equation}
where $\omega_h(\xi)\in {\tilde S}^{\rho}$ and
\[
\lr{\xi}^{\rho}_{\ell}/4\leq \omega_h(\xi)\leq \lr{\xi}_{\ell}^{\rho}.
\]
\end{lemma}
{\bf Proof.} Denote 
\[
\frac{W_h^{n+1}-W^n_h}{k}=-2\,{\tilde \omega}_h \frac{W^{n+1}_h+W^n_h}{2}
\]
then it is clear that
\[
{\tilde \omega}_h=\frac{1-e^{-ak\lr{\xi}^{\rho}_{\ell}\chi_h}}{k}\,\frac{1}{1+e^{-ak\lr{\xi}^{\rho}_{\ell}\chi_h}}.
\]
Since
\[
\frac{1-e^{-ak\lr{\xi}^{\rho}_{\ell}\chi_h}}{k}=a\lr{\xi}_{\ell}^{\rho}\chi_h\int_0^1e^{-ak\theta\lr{\xi}^{\rho}_{\ell}\chi_h}d\theta
\]
one can define $\omega_h$ by
\[
{\tilde \omega}_h=a\Big(\lr{\xi}_{\ell}^{\rho}\int_0^1e^{-ak\theta\lr{\xi}^{\rho}_{\ell}\chi_h}d\theta\,\frac{1}{1+e^{-ak\lr{\xi}^{\rho}_{\ell}\chi_h}}\Big)\chi_h=a\,\omega_h\,\chi_h.
\]
Then the first assertion is clear from \eqref{eq:W_h}.
 Note that 
 \begin{equation}
 \label{eq:estW}
 \big|\partial_{\xi}^{\alpha}\big(a\,k\lr{\xi}_{\ell}^{\rho}\chi_h \big)\big|\leq C_{\alpha}\lr{\xi}_{\ell}^{-|\alpha|}
 \end{equation}
  because of \eqref{eq:constraint:1}. Therefore   one has $
 |\partial_{\xi}^{\alpha}\omega_h|\leq C_{\alpha}\, \lr{\xi}^{\rho-|\alpha|}_{\ell}
$.
Using  \eqref{eq:W_h} 
this implies the second assertion.
\hfill
\qed
%

\subsection{Crank-Nicholson after conjugation}

Note that $u^n_h$ satisfy 
\begin{equation}
\label{eq:h2h}
\forall\, n\in\NN,\qquad \chi_hu^n_h=u^n_h
\end{equation}
thanks to  \eqref{eq:smallspec}. Assume that $u^n_h$ satisfies
\begin{equation}
\label{eq:eqCr_Ni}
\delta_ku^n_h=\frac{u^{n+1}_h-u^n_h}{k}=G^n_h(x,D)\ \frac{u^{n+1}_h+u^n_h}{2}+f^n
\end{equation}
where $\chi_hf^n=f^n$ is not necessarily assumed.

Consider a weighted energy $(R_h^nW^n_hu_h^n\,,\,W^n_hu_h^n)$ where $R^n_h$ is a symmetrizer that  is symmetric $(R^n_h)^*=R^n_h$
 and will be defined in Subsection \ref{sec:defR} below. 
The discrete analog of $\partial_t (R_h^nW^n_hu_h^n\,,\,W^n_hu_h^n)$ is 
the time difference 
\begin{equation}
\label{eq:timediff}
\begin{split}
\delta_k&(R^n_hW^n_h\,u^n_h,\,W^n_h\,u^n_h) \\ 
&=
\frac{
(R^{n+1}_hW^{n+1}_h
u^{n+1}_h,\,W^{n+1}_hu^{n+1}_h)-(R^{n}_hW^n_h
u^{n}_h,\,W^n_hu^{n}_h)}{k}.
\end{split}
\end{equation}

In what follows we omit the subscript $h$ for 
ease of reading.  Write \eqref{eq:timediff} as
\begin{align*}
\frac{
(W^{n+1}R^{n}W^{n+1}
u^{n+1}
,\, u^{n+1})
-
(W^nR^{n}W^n
u^{n}
,\, u^{n})
}
{k}
\ +\ (III)\,,
\end{align*}
with 
\begin{equation}
\label{eq:firsterror}
(III)\ :=\ \frac{
((R^{n+1}-R^n)W^{n+1}
u^{n+1},\,W^{n+1}u^{n+1})}{k}.
\end{equation}
The term $(III)$ 
is an error term that will be estimated in Subsection \ref{sec:defR}. The first term  is 
   equal to
\begin{equation}
\label{eq:3lines}
\begin{split}
&\frac{
((R^n{\bar\delta}_k{W^n}) u^{n+1},\,
u^{n+1})
+
( (R^n{\bar\delta}_k{W^n})u^{n},\,
u^{n})
}
{2}
\\
&+\ \bigg(\Big(\frac{W^{n+1}R^{n}W^{n+1}+W^nR^{n}W^n}{2}\Big)
 \Big(\frac{u^{n+1}+u^{n}}{2}\Big)\ ,\ 
\delta_ku^n
\bigg)
\\
&+\ \bigg(
\Big(\frac{W^{n+1}R^{n}W^{n+1}+W^nR^{n}W^n}{2}\Big)\
\delta_ku^n\,,\,
\Big(\frac{u^{n+1}+u^{n}}{2}\Big)
\bigg)
\end{split}
\end{equation}
where
\[
R^n{\bar\delta}_k{W^n}=\frac{W^{n+1}R^nW^{n+1}-W^nR^nW^n}{k}.
\]
The first line of \eqref{eq:3lines} is
\begin{align}
(I)=\frac{
((R^n{\bar\delta}_k{W^n})u^{n+1},
u^{n+1})
+
( (R^n{\bar\delta}_k{W^n})u^n,
u^{n})
}
{2}\,.
\end{align}
Note that
\begin{equation*}
\label{eq:E12}
\begin{split}
R^n{\bar\delta}_k{W^n}
=\frac{W^{n+1}R^nW^{n+1}-W^nR^nW^n}{k}
\\
=\frac{1}{2}\,\frac{W^{n+1}-W^n}{k}R^n(W^{n+1}+W^n)
+\frac{1}{2}\,(W^{n+1}+W^n)R^n\frac{W^{n+1}-W^n}{k}.
\end{split}
\end{equation*}

Using \eqref{eq:W_CrNi} and $\omega\chi_h\,W^m=W^m \omega\chi_h$
this becomes
\begin{align*}
R^n{\bar\delta}_k{W^n}
=
-\frac{a}{2}\,(W^{n+1} & +W^n)\,\omega\,\chi_hR^n (W^{n+1}+W^n)
\\
&-\frac{a}{2}\,(W^{n+1}+W^n)R^n\, \omega\,\chi_h\,(W^{n+1}+W^n) .
\end{align*}
Therefore with $\Omega^{n}=W^{n+1}+W^n$ one has, since   $(R^n)^*=R^n$
\begin{align*}
( (R^n{\bar\delta}_k{W^n})w,
w)=-a\,{\mathsf{Re}}\, ( R^n\, \Omega^n\, w, \,\omega\chi_h\, \Omega^n w)\\
=-a\,{\mathsf{Re}}\,(\omega\chi_h\,  R^n\, \Omega^n\, w, \, \Omega^n w).
\end{align*}
Thus (I) yields
\begin{align*}
(I)\ =\
-a\sum_{j=0}^1\,{\mathsf{Re}}\, (\omega\chi_h\,R^n\,\, \Omega^n\, u^{n+j},\Omega^n\, u^{n+j}).
\end{align*}
Since $\chi_h\,\Omega^n=\Omega^n\,\chi_h$ and $\omega\chi_h=\chi_h \omega$ 
and using $\chi_hu^{n+j}=u^{n+j}$ that follows from \eqref{eq:h2h} one has
\begin{equation}
\label{eq:mainCrNi}
(I)\ =\ 
-a\,\sum_{j=0}^1\,{\mathsf{Re}}\, (\omega\, R^n\, \Omega^n u^{n+j},\,\Omega^n u^{n+j}).
\end{equation}
 The second line of \eqref{eq:3lines} yields, with $U^n=u^{n+1}+u^n$ 
\begin{align*}
\bigg(
\Big(\frac{W^{n+1}R^nW^{n+1}+W^n R^nW^n}{2}\Big)\Big(\frac{u^{n+1}+u^{n}}{2}\Big), \delta_ku^n
\bigg)\\
=\frac{1}{8}(R^nW^{n}U^n,\, W^{n}G^nU^n)+\frac{1}{8}(R^nW^{n+1}U^n,\, W^{n+1}G^nU^n)\\
+\frac{1}{4}(U^n,\, (W^{n+1}R^nW^{n+1}+W^n R^nW^n) f^n)\,.
\end{align*}
Because of \eqref{eq:eqCr_Ni}, this is equal to
\begin{align*}
\frac{1}{8}\sum_{j=0}^1( R^n\,W^{n+j} \,U^n,\, W^{n+j}\,G^n\,U^n)
+\frac{1}{4}\sum_{j=0}^1(U^n, \,W^{n+j}R^n\,W^{n+j}f^n).
\end{align*}
Similarly  the third line of \eqref{eq:3lines} is
\begin{align*}
\Big(\frac{W^{n+1}R^nW^{n+1}+W^nR^{n}W^n}{2}\ 
\delta_ku^n,\,
\frac{u^{n+1}+u^{n}}{2}
\Big)&
\\
=\frac{1}{8}\sum_{j=0}^1(W^{n+j}\,G^n\,U^n,\, R^n\,W^{n+j}\,U^n)
+\frac{1}{4}\sum_{j=0}^1  &   (W^{n+j}R^nW^{n+j} f^n,\, U^n).
\end{align*}
Therefore the sum of the second and the third lines of \eqref{eq:3lines}, denoted by $(II)$,  yields
\begin{equation}
\label{eq:subCrNi}
\begin{split}
(II)=\frac{1}{4}\sum_{j=0}^1\,{\mathsf{Re}}\,( R^n & \,W^{n+j} \,U^n,\, W^{n+j}\,G^n\,U^n)\\
&+\frac{1}{2}\sum_{j=0}^1\,{\mathsf{Re}}\,(U^n, \,W^{n+j}R^n\,W^{n+j}f^n).
\end{split}
\end{equation}
Recalling 
\begin{equation}
\label{eq:I+II+III}
\delta_k(R^nW^n\,u^n,\,W^n\,u^n)=(I)+(II)+(III)
\end{equation}
we have proved the following proposition.
\begin{proposition}
\label{pro:I+II+III}We have
\begin{align*}
\delta_k(R^nW^n\,u^n,\,W^n\,u^n)=-a\,\sum_{j=0}^1\,{\mathsf{Re}}\, (\omega\, R^n\, \Omega^n u^{n+j},\,\Omega^n u^{n+j})\\
+\frac{1}{4}\sum_{j=0}^1\,{\mathsf{Re}}\,( R^n\,W^{n+j} \,U^n,\, W^{n+j}\,G^n\,U^n)
+\frac{1}{2}\sum_{j=0}^1\,{\mathsf{Re}}\,(U^n, \,W^{n+j}R^n\,W^{n+j}f^n)\\
+\frac{
((R^{n+1}-R^n)W^{n+1}
u^{n+1},\,W^{n+1}u^{n+1})}{k}
\end{align*}
where $\Omega^n:=W^{n+1}+W^n$ and $U^n:=u^{n+1}+u^n$.
\end{proposition}
%

\subsection{Composition with $W^{\pm n}_h $ and definition of $R^n_h$} 
\label{sec:defR}

First recall 
Definition 2.4 from
\cite{CNR}.

\begin{definition}
\label{dfn:kurasu}\rm
For  $1<s$, $m\in\RR$, 
the family $a(x,\xi;\ell)\in C^{\infty}(\RR^d\times\RR^d ) $ 
belongs to  ${\tilde S}_{(s)}^m$ if there exist $C>0$, $A>0$ independent of 
$\ell\ge1,x,\xi$ such that for all $\alpha$, $\beta\in\NN^d$,
\[
\big|
\partial_x^{\beta}\partial_{\xi}^{\alpha}a(x,\xi;\ell)
\big|
\ \leq\ 
C\,A^{|\alpha+\beta|}\ 
|\alpha+\beta|!^s\
\lr{\xi}_{\ell}^{m-|\alpha|}
\,.
\]
\end{definition}

We often write $a(x,\xi)$ for $a(x,\xi;\ell)$ dropping the $\ell$. If $a(x,\xi)$ is the symbol of a differential operator of order $m$ with coefficients $a_{\alpha}(x)\in G^{s}(\RR^d)$ 
then $a(x,\xi)\in {\tilde S}^m_{(s)}$ because $|\partial_{\xi}^{\beta}\xi^{\alpha}|\leq CA^{|\beta|}|\beta|!\lr{\xi}_{\ell}^{|\alpha|-|\beta|}$ and $|\partial_x^{\beta}a_{\alpha}(x)|\leq C_{\alpha}A_{\alpha}^{|\beta|}|\beta|!^s$ for any $\beta\in\NN^d$. 

\smallskip

\begin{proposition}
\label{pro:b_1}Suppose $1/2\leq \rho<1$ and $s=1/\rho$ and let $A(x,\xi)$ be $m\times m$ matrix valued with entries in ${\tilde S}^1_{(s)}$ and $\partial_x^{\al}A(x,\xi)=0$ outside $|x|<R$ for some $R>0$ if $|\al|>0$. 
Define $
m^*:={\max{\{\rho-k(1-\rho),-1+\rho\}}}$.
Then there is ${\bar\tau}>0$,  $\ell_0>0$ such that  
\[
{\tilde A}(x,D)=e^{\tau\lr{D}_{\ell}^{\rho}\chi_h}A(x,D)e^{-\tau\lr{D}_{\ell}^{\rho}\chi_h}
\]
 is a pseudodifferential operator with symbol given by
\[
{\tilde A}(x,\xi)=\sum_{|\alpha|\leq k}\frac{1}{\al!}D_x^{\alpha}A(x,\xi)\big(\tau\,\nabla_{\xi}(\lr{\xi}_{\ell}^{\rho}\chi_h)\big)^{\alpha}+R_k(x,\xi)
\]
with $R_k\in {\tilde S}^{m^*}$ uniformly in $\tau$, $\ell$  constrained to satisfy 
\begin{equation}
\label{eq:constraint:3}
|\tau|\leq  {\bar\tau},\qquad \ell\geq \ell_0\,.
\end{equation}
 In particular ${\tilde A}(x,\xi)\in {\tilde S}^1$ uniformly in such
  $\tau$, $\ell$.
 \end{proposition}
\begin{remark}
\label{rem:spcut}\rm This proposition  with $\chi_h\equiv 1$
 is \cite[Proposition 2.6]{CNR}. 
 The proof for the case  $\chi_h\equiv 1$  works 
 without any change for the case  
 $\chi_h\in {\tilde S}^0$.
\end{remark}
Choosing a smaller ${\bar\tau}>0$ if necessary one can assume that
\[
{\bar\tau}\big|\nabla_{\xi}(\lr{\xi}_{\ell}^{\rho}\chi_h)\big|
\ \leq\ 
 \lr{\xi}_{\ell}^{\rho-1}.
\]

In what follows we choose $T={\bar \tau}$ in the definition of $W(t,\xi)$ 
yielding
\[
W(t,\xi)
\ =\
e^{({\bar\tau}-t)\lr{\xi}^{\rho}_{\ell}\chi_h}.
\]
 With $N=\max\{2\theta,m\}$
 denote
\[
H( t,x,\xi, \tau)
\ =\
\sum_{|\alpha|\leq  N}\frac{1}{\alpha !}
D_x^{\alpha}A(t,x,\xi)
\big(({\bar\tau}-\tau)\nabla_{\xi}(\lr{\xi}_{\ell}^{\rho}\chi_h)\big)^{\alpha}
\,.
\]
Suppressing the subscript $h$ for ease of reading, Proposition \ref{pro:b_1} shows
that
\[
W(\tau,\xi)\#A(t,x,\xi)\#W^{-1}(\tau,\xi)
\ =\
H(t,x,\xi,\tau)+R,\quad R\in {\tilde S}^{m^*} \,.
\]
The choice of $N$ guarantees that 
where $2\theta(1-\rho)+m^*\leq \rho$.
 Define 
 \[
 H^h(t,x,\xi,\tau)
 \ =\
 \chi_{2h}^2(\xi)H(t,x,\xi,\tau)\,.
 \]
Then, the definition of ${\mathcal H}_N$ implies that 
\begin{equation}
\label{eq:defHh}
\begin{split}
&H^h(t,x,\xi,\tau)\\
&=\chi_{2h}^2(\xi)\,
\lr{\xi}_\ell\ {\mathcal H}_N\big( t,x,\xi/\lr{\xi}_\ell, ({\bar\tau}-\tau)\nabla_{\xi}(\lr{\xi}_{\ell}^{\rho}\chi_h)/\lr{\xi}^{\rho-1}_{\ell}, 0, \lr{\xi}_\ell^{\rho-1} \big).
\end{split}
\end{equation}
 In \eqref{eq:defHh}
choose
\begin{equation}
\label{eq:choice}
\begin{split}
s=\chi_{2h}^2(\xi)\lr{\xi}_\ell, \quad \epsilon=\lr{\xi}_\ell^{\rho-1}, \quad \xi=\xi/\lr{\xi}_\ell, \\
y=({\bar\tau}-\tau)\nabla_{\xi}(\lr{\xi}_{\ell}^{\rho}\chi_h)/\lr{\xi}^{\rho-1}_{\ell}, \quad \eta=0\,.
\end{split}
\end{equation}
 Using $0\leq \chi_{2h}\leq 1$, it follows from \eqref{eq:Matexp} that
\begin{equation}
\label{eq:expH}
\begin{split}
 \lr{\xi}_\ell^{-\theta(1-\rho)}  e^{-cs\lr{\xi}_\ell^{\rho}}/C
\ \leq\ 
\big\| e^{isH^h(t,x,\xi,\tau)} \big\|
 \leq
C \lr{\xi}_\ell^{\theta(1-\rho)}    e^{cs\lr{\xi}_\ell^{\rho}}  
\end{split}
\end{equation}
for $|t|\leq T$,  $\ell\geq \ell_0$ where 
\begin{equation}
\label{eq:constraint6}
0\leq \tau\leq {\bar\tau},\quad 0< \epsilon=\lr{\xi}_{\ell}^{\rho-1}
 \leq \ell^{\rho-1}\leq \ell_0^{\rho-1}=\epsilon_0.
\end{equation}
Following \cite{CNR}  define
\[
M^h(t,x,\xi,\tau)
\ =\
iH^h(t,x,\xi,\tau)
\ -\
b\, \lr{\xi}_\ell^{\rho}
\]
with a positive parameter $b>0$ that will be fixed later. 
Since $\|e^{sM^h(t,x,\xi,\tau)}\|=e^{-bs\lr{\xi}_{\ell}^{\rho}}\|e^{isH^h(t,x,\xi,\tau)}\|$,  \eqref{eq:expH} implies
\[
\lr{\xi}_{\ell}^{-\nu}\,e^{-c_1b\,s\lr{\xi}_{\ell}^{\rho}}/C\leq \|e^{sM^h(t,x,\xi,\tau)}\|\leq C\,\lr{\xi}_{\ell}^{\nu}\,e^{-c_2b\,s\lr{\xi}_{\ell}^{\rho}}
\]
with $\nu=\theta(1-\rho)$ for $|t|\leq T$ and $b\geq b_0$ with some $b_0>0$ where  $c_1$, $c_2$ and $C>0$ are independent of $\ell$, $h$ and $b$.

Introduce the symmetrizer
$$
R_h(t,x,\xi,\tau)
\ :=\
b\int_0^\infty
\langle\xi\rangle^\rho_\ell
\big(
e^{sM^h(t,x,\xi,\tau)}\big)^*
\big(
e^{sM^h(t,x,\xi,\tau)}\big)
ds
\,.
$$
From \cite[Theorem 3.1]{CNR} it follows that
\[
R_h(t,x,\xi,\tau)\in {\tilde S}^{2\nu}_{\rho-\nu,1-\rho+\nu},\qquad b\,\partial_tR_h(t,x,\xi,\tau)\in {\tilde S}^{3\nu+1-\rho}_{\rho-\nu,1-\rho+\nu}
\]
under the constraint
\begin{equation}
\label{eq:constraint:b}
 b\,\ell^{-(1-\rho)}\leq 1
\end{equation}
so that $b\,\lr{\xi}_{\ell}^{\rho-1}\leq b\,\ell^{-(1-\rho)}\leq 1$. 
Recall \cite[page 230]{CNR} that
\begin{align*}
R_hM^h+(M^h)^*R_h
=R_h\big(i\chi_{2h}^2H-b\lr{\xi}^{\rho}_{\ell}\big)+\big(i\chi_{2h}^2H-b\lr{\xi}^{\rho}_{\ell}\big)^*R_h=-b\,\lr{\xi}^{\rho}_{\ell}
\end{align*}
that is
\begin{equation}
\label{eq:RiH}
R_h(i\chi_{2h}^2H)+(i\chi_{2h}^2H)^*R_h=-b\,\lr{\xi}_{\ell}^{\rho}+2b\,\lr{\xi}_{\ell}^{\rho}\,R_h.
\end{equation}
\begin{lemma}
\label{lem:Rn+1Rn:2}We have
\[
\frac{b\,\big(R_h(t,x,\xi,a(n+1)k)-R_h(t,x,\xi,ank)\big)}{k\,a}\in {\tilde S}^{3\nu}_{\rho-\nu,1-\rho+\nu}
\]
for $0\leq t\leq T$ uniformly in $a, b$, $n, k, h$ under the constraint $
 ank\leq {\bar\tau}$.
\end{lemma}
{\bf Proof.} We show that
\begin{equation}
\label{eq:diff_nR}
\big|\partial_x^{\beta}\partial_{\xi}^{\alpha}\partial_{\tau}R_h(t,x,\xi,\tau)\big|\leq C_{\alpha\beta}\,b^{-|\alpha+\beta|-1}\lr{\xi}_{\ell}^{3\nu+(1-\rho+\nu)|\beta|-(\rho-\nu)|\alpha|}
\end{equation}
with $C_{\alpha \beta}$ independent of $b$, $ h$ and $0\leq \tau\leq {\bar\tau}$. 
If \eqref{eq:diff_nR} is proved then writing
\[
R_h(t,x,\xi,a(n+1)k)-R_h(t,x,\xi,ank)=\int_{ank}^{ank+ak}\partial_{\tau}R_h(t,x,\xi,\nu)d\nu
\]
the assertion follows immediately. To prove the estimate \eqref{eq:diff_nR} we apply the same arguments in the proof of \cite[Theorem 3.1]{CNR}.
First consider $\partial_{\tau}H(t,x,\xi,\tau)$. 
Since
\[
\partial_{\tau}H^h(t,x,\xi,\tau)=-\chi_{2h}^2\sum_{1\leq |\alpha|\leq  N}({\bar\tau}-\tau)^{|\alpha|-1}\frac{|\al|}{\alpha !}
D_x^{\alpha}A(t,x,\xi)
\big(\nabla_{\xi}(\lr{\xi}_{\ell}^{\rho}\chi_h)\big)^{\alpha}
\]
it follows that
\begin{equation}
\label{eq:diff_nH}
\big|\partial_x^{\beta}\partial_{\xi}^{\alpha}\partial_{\tau}H^h(t,x,\xi,\tau)\big|\leq C_{\alpha\beta}\lr{\xi}_{\ell}^{\rho-|\alpha|}.
\end{equation}
Denote 
\[
X(s;t,x,\xi,\tau)=e^{sM^h(t,x,\xi,\tau)}v,\;\;v\in\CC^m,\quad X^{\alpha}_{\tau\beta}=\partial_x^{\beta}\partial_{\xi}^{\alpha}\partial_{\tau}X(t,x,\xi,\tau).
\]
Since 
\[
{\dot X}_{\tau}=M^hX_{\tau}+\partial_{\tau}H^hX,\quad X_{\tau}(0)=0
\]
then \eqref{eq:diff_nH} and Duhamel's representation yields
\[
\big|X_{\tau}\big|=\Big|\int_0^se^{(s-{\tilde s})M^h}(\partial_{\tau}H^h)Xd{\tilde s}\Big|\leq C(s+\lr{\xi}_{\ell})\lr{\xi}_{\ell}^{\nu+\rho-1}E(s)
\]
where $E(s)=\lr{\xi}_{\ell}^{\nu}e^{-cb s\lr{\xi}_{\ell}^{\rho}}$. Repeating the same arguments in the proof of \cite[Theorem3.1]{CNR} one can prove
\[
\big|X^{\alpha}_{\tau\beta}\big|\leq C_{\alpha\beta}\,(s+\lr{\xi}_{\ell}^{-1})^{|\alpha|}(1+s\lr{\xi}_{\ell})^{|\beta|+1}\lr{\xi}_{\ell}^{\nu(|\alpha+\beta|+1)+\rho-1}E(s)
\]
from which we obtain \eqref{eq:diff_nR} by exactly the same way as in the proof of \cite[Theorem 3.1]{CNR}. 
\hfill
\qed
\vskip.2cm
\begin{lemma}
\label{lem:Rn+1Rn}
With 
$
R_h^n(x,\xi)
\, :=\,
R_h(nk,x,\xi,ank)
$\,,
one has
\[
\frac{b\,\big(R_h^{n+1}(x,\xi)-R_h^n(x,\xi)\big)}{k}\in {\tilde S}^{-2\nu+\rho}_{\rho-\nu,1-\rho+\nu}
\]
for $0\leq (n+1)k\leq T$ uniformly in $a, b$, $n, k, h$ under the constraint 
\begin{equation}
\label{eq:constraint:a}
ank\leq {\bar\tau},\qquad a\,\ell^{-\rho/6}\leq 1.
\end{equation}
\end{lemma}
{\bf Proof.} Write
\begin{align*}
R^{n+1}_h-R^n_h=R_h((n+1)k,x,\xi,a(n+1)k)-R_h((n+1)k,x,\xi,ank)\\
+R_h((n+1)k,x,\xi,ank)-R_h(nk,x,\xi,ank).
\end{align*}
Express
\begin{align*}
R_h((n+1)k,x,\xi,ank)-R_h(nk,x,\xi,ank)
=\int_{nk}^{nk+k}\partial_tR_h(t',x,\xi,ank)dt'
\,.
\end{align*}
Using
 $b\,\partial_tR_h(t,x,\xi,\tau)\in {\tilde S}^{3\nu+1-\rho}_{\rho-\nu,1-\rho+\nu}$, one obtains 
\[
\frac{b\,(R_h((n+1)k,x,\xi,ank)-R_h(nk,x,\xi,ank))}{k}
\ \in\ {\tilde S}^{3\nu+1-\rho}_{\rho-\nu,1-\rho+\nu}
\]
where $3\nu+1-\rho\leq -2\nu+\rho$ in view of  \eqref{rem:rhoandnu}. For the term $R_h((n+1)k,x,\xi,a(n+1)k)-R_h((n+1)k,x,\xi,ank)$ we apply Lemma \ref{lem:Rn+1Rn:2} to get 
\[
\frac{b\,\big(R_h((n+1)k,x,\xi,a(n+1)k)-R_h((n+1)k,x,\xi,ank)\big)}{k\,a}\in {\tilde S}^{3\nu}_{\rho-\nu,1-\rho+\nu}.
\]
Here note that \eqref{eq:rhoandnu} implies $
1/2>3\nu$ 
because $1>\rho\geq 3\nu+1/2$ and hence 
\begin{equation}
\label{eq:rhonu}
\rho\geq 3\nu+1/2>6\nu.
\end{equation}
Then noting that $a\lr{\xi}_{\ell}^{3\nu}\leq a\lr{\xi}^{-\rho/6}_{\ell}\lr{\xi}^{\rho-2\nu}_{\ell}\leq a\,\ell^{-\rho/6}\lr{\xi}^{\rho-2\nu}_{\ell}$ one has
\[
\frac{b\,\big(R_h((n+1)k,x,\xi,a(n+1)k)-R_h((n+1)k,x,\xi,ank)\big)}{k}\in {\tilde S}^{-2\nu+\rho}_{\rho-\nu,1-\rho+\nu}
\]
under the constraint $a\,\ell^{-\rho/6}\leq 1$. Thus the proof is complete.
\hfill\qed\vskip.2cm

\begin{definition}
\label{dfn:pos_definite}\rm 
For a $m\times m$ complex matrix $\caM=\caM^*$,  the notation $\caM\gg 0$ means that for all
$v\in \CC^m$ one has $(\caM v\,,\,v)_{\CC^m}\ge 0$.
For two such matrices,  $\caM_1\gg\caM_2$ means that $\caM_1-\caM_2 \gg 0$.
\end{definition}
%
Equation
\eqref{eq:expH} yields for any $v\in \CC^m$
\begin{align*}
&(R^n_h(x,\xi)v,v) =
b  \int_0^{\infty}\lr{\xi}_{\ell}^{\rho}\|e^{sM^h(nk,x,\xi,ank)}v\|^2\,ds\\
 &\geq\  C^{-2}\,\|v\|^2\lr{\xi}_{\ell}^{-2\nu}
 \int_0^{\infty}b\,\lr{\xi}_\ell^{\rho}\,e^{-2c_1b\,s\lr{\xi}_\ell^{\rho}}\,ds
 \geq
 c\,
\lr{\xi}_\ell^{-2\nu}\|v\|^2\,.
\end{align*}
This is an important pointwise lower bound for the symbol
\begin{equation}
\label{eq:lowerboundR}
R^n_h(x,\xi)
\ \gg\
 c\,\lr{\xi}_{\ell}^{-2\nu}I
\end{equation}
where $c$ is independent of $b$, $a$, $n$, $k$, $h$ constrained to satisfy  \eqref{eq:constraint:a} and \eqref{eq:constraint:b}.

\subsection{Estimate of (I)}

Suppressing the suffix $h$ again,
denote
\[
{ W}(\xi)=e^{-ak\lr{\xi}_{\ell}^{\rho}\chi_h}
\]
so that $W^{n+1}=W^nW$ where $1/2\leq W\leq 1$ and $W^{\pm 1}\in {\tilde S}^0$ which follows from  \eqref{eq:W_h} and \eqref{eq:estW}. Consider ${\mathsf{Re}}(\omega\,R^n\,\Omega^n\,w,\, \Omega^n\,w)$. Write $\Omega^n=W^{n+1}+W^n=(1+W)\,W^n$ and hence
\[
{\mathsf{Re}}\,(\omega\,R^n\,\Omega^n\,w,\Omega^n\,w)
\ =\
{\mathsf{Re}}\,\big((1+W)\omega\,R^n(1+W)W^n\,w,\,W^n\,w\big)\,.
\]
Note that $R^n\in {\tilde S}^{2\nu}_{\rho-\nu,1-\rho+\nu}$ uniformly for  parameters 
satisfying the constraints \eqref{eq:constraint:a} and \eqref{eq:constraint:b}.  
\begin{lemma}
\label{lem:remcompo}One can write
\[
(1+W)\#\omega\#R^n\#(1+W)=(1+W)^2\,\omega\,R^n+iR_1^n+R_2^n
\]
with $(R_1^n)^*=R_1^n$ and $R_2^n\in {\tilde S}^{4\nu-\rho}_{\rho-\nu,1-\rho+\nu}$.
\end{lemma}
\noindent
{\bf Proof.} Denote $f(\xi)=(1+W(\xi))\omega(\xi)$ and $g(\xi)=1+W(\xi)$. Since one has $f(\xi)\in {\tilde S}^{\rho}\subset {\tilde S}^{\rho}_{\rho-\nu,1-\rho+\nu}$ and $g(\xi)\in {\tilde S}^0\subset {\tilde S}^{0}_{\rho-\nu,1-\rho+\nu}$ applying \cite[Theorem 18.5.4]{Hobook} one can write
\begin{align*}
\big((1+W)\omega\big)\#R^n\#(1+W)=(1+W)^2\omega R^n\\
+\frac{1}{2i}\sum_{|\alpha+\beta|=1}(-1)^{|\beta|}\partial_{\xi}^{\alpha}f\,(\partial_x^{\alpha+\beta}R^n)\,\partial_{\xi}^{\beta}g\\
+\sum_{2\leq |\alpha+\beta|<N}\frac{(-1)^{|\beta|}}{(2i)^{|\alpha+\beta|}\alpha!\beta!}\partial_{\xi}^{\alpha}f\,(\partial_x^{\alpha+\beta}R^n)\,\partial_{\xi}^{\beta}g+R_N
\end{align*}
where $R_N\in {\tilde S}^{2\nu+\rho-N(2\rho-1-2\nu)}_{\rho-\nu,1-\rho+\nu}$. 

Choose $N$ so large that
 $2\nu+\rho-N(2\rho-1-2\nu)\leq 4\nu-\rho$. The second term on the right-hand side, denoted by $iR_1^n$, clearly satisfies $(R_1^n)^*=R_1^n$ because $f(\xi)$ and $g(\xi)$ are real scalar symbols. Since $\partial_{\xi}^{\alpha}f\in {\tilde S}^{\rho-|\alpha|}_{\rho-\nu,1-\rho+\nu}$ and $\partial_{\xi}^{\beta}g\in {\tilde S}^{-|\beta|}_{\rho-\nu,1-\rho+\nu}$ it is clear that the third term on the right-hand side is in ${\tilde S}^{4\nu-\rho}_{\rho-\nu,1-\rho+\nu}$. 
\hfill
\qed
\vskip.2cm

Thanks to Lemma \ref{lem:remcompo} we have
\begin{align*}
{\mathsf{Re}}((1+& W)\,  \omega R^n(1+W)W^n\,w,\,W^n\,w)
\\
&\geq\
 \big({\rm Op}((1+W)^2\omega R^n)\,W^n\,w,\,W^n\,w\big)
-C\|\lr{D}_{\ell}^{2\nu-\rho/2}\,W^n\,w\|^2.
\end{align*}

It follows from Lemma \ref{lem:hyoome} and \eqref{eq:lowerboundR} that
\[
(1+W)^2\omega R^n\gg c \big(\lr{\xi}_{\ell}^{\rho-2\nu}I+\lr{\xi}^{\rho}_{\ell}R^n\big),\quad (1+W)^2\omega R^n\in {\tilde S}^{2\nu+\rho}_{\rho-\nu,1-\rho+\nu}
\]
 with  $c>0$ uniformly in the constrained parameters
 $k, n, h$, $a, b$.  Note that 
\[
(1+W)^2\omega R^n
\ \in\
 S\big(\lr{\xi}_{\ell}^{2\nu+\rho},\, b^{-2}(\lr{\xi}_{\ell}^{2(1-\rho+\nu)}|dx|^2+\lr{\xi}_{\ell}^{-2(\rho-\nu)}|d\xi|^2\big)
\]
since for any $k(\xi) \in {\tilde S}^{q}$ one has 
\[
|\partial_{\xi}^{\alpha}k(\xi)|\leq C_{\alpha}\lr{\xi}_{\ell}^{q-|\alpha|}\leq C_{\alpha}\ell^{-|\alpha|(1-\rho+\nu)}\lr{\xi}_{\ell}^{q-|\alpha|(\rho-\nu)}
\]
 which is bounded by $C_{\alpha}b^{-|\alpha|}\lr{\xi}_{\ell}^{q-|\alpha|(\rho-\nu)}$ because of  \eqref{eq:constraint:b}. 
Repeating the  arguments proving \cite[(4.6)]{CNR} it follows from the sharp G\aa rding inequality \cite[Theorem 18.6.7]{Hobook} that there is $\ell_0>0$ such that
\begin{align*}
({\rm Op} & ((1+W)^2\omega R^n)\,W^n\,w,W^n\,w)\geq c\,\|\lr{D}_{\ell}^{-\nu+\rho/2}W^n\,w\|^2
\\
&+c\,({\rm Op}(\lr{\xi}_{\ell}^{\rho}R^n)\,W^n\,w,W^n\,w)-Cb^{-2}\|\lr{D}_{\ell}^{-\rho/2+1/2+2\nu}W^n\,w\|^2
\end{align*}
for  $\ell\geq \ell_0$. 
Since $-\nu+\rho/2\geq -\rho/2+1/2+2\nu$ by \eqref{eq:rhoandnu},
 choosing another $b_0$ if necessary one obtains
\begin{equation}
\label{eq:rhoRn}
\begin{split}
({\rm Op}((1 & +W)^2\omega R^n)\,W^n\,w,\,W^n\,w)
\\
&\geq c'\,\|\lr{D}_{\ell}^{-\nu+\rho/2}W^n\,w\|^2
+c'\,({\rm Op}(\lr{\xi}^{\rho}_{\ell}R^n)\,W^n\,w,\,W^n\,w)
\end{split}
\end{equation}
with $c'>0$ uniform for $\ell\geq \ell_0$ and $b\geq b_0$. From \eqref{eq:rhoandnu} again  one sees that 
\[
\|\lr{D}_{\ell}^{2\nu-\rho/2}\,v\|^2
\ \leq\
 C\ell^{-1}\|\lr{D}_{\ell}^{-\nu+\rho/2}\,v\|^2
\]
and one concludes that choosing another $\ell_0$ if necessary 
\begin{align*}
{\mathsf{Re}}\,((1&+W)\omega\,R^n\,(1+W)W^n\,w,\,W^n\,w)
\\
&\geq c\,\|\lr{D}_{\ell}^{-\nu+\rho/2}W^n\,w\|^2
+c\,({\rm Op}(\lr{\xi}^{\rho}_{\ell}R^n)\,W^n\,w,\,W^n\,w)
\end{align*}
with  $c>0$ uniform for $\ell\geq \ell_0$ and $b\geq b_0$. Repeating the same arguments one obtains
\begin{align*}
{\mathsf{Re}}\,((1&+W^{-1})\,\omega\,R^n
(1+W^{-1})\,W^{n+1}\,w,\,W^{n+1}w)
\\
&\geq c\,\|\lr{D}_{\ell}^{-\nu+\rho/2}W^{n+1}\,w\|^2
+c\,({\rm Op}(\lr{\xi}^{\rho}_{\ell}R^n)\,W^{n+1}\,w,\,W^{n+1}\,w).
\end{align*}
Summarizing we have
\begin{lemma}
\label{lem:yoiko}There are $c>0$, $\ell_0$ and $b_0$  such that for $\ell\geq \ell_0$ and $b\geq b_0$ one has
\begin{align*}
{\mathsf{Re}}\,( & \omega\,R^n\,\Omega^n\,w,\,\Omega^n\,w)
\\
& \geq c\,\Big(\sum_{i=0}^1\|\lr{D}_{\ell}^{-\nu+\rho/2}\,W^{n+i}\,w\|^2
+\big({\rm Op}(\lr{\xi}^{\rho}_{\ell}R^n)\,W^{n+i}\,w,\,W^{n+i}\,w\big)\Big).
\end{align*}
\end{lemma}
 Lemma \ref{lem:yoiko} together with \eqref{eq:mainCrNi} prove the following
\begin{proposition}
\label{pro:iiko}There exist $c>0$, $\ell_0>0$ and $b_0$  such that for $\ell\geq \ell_0$ and $b\geq b_0$ one has
\begin{align*}
(I)  \leq 
-c\,a &
\sum_{0\le i\le1\atop 0\le j\le 1}
\Big(
\|\lr{D}_{\ell}^{-\nu+\rho/2}\,W^{n+i}\,u^{n+j}\|^2
\\ 
\ &+\ 
\big({\rm Op}(\lr{\xi}^{\rho}_{\ell}R^n)\,W^{n+i}\,u^{n+j},\,W^{n+i}\,u^{n+j}\big)
\Big)\,.
\end{align*}
\end{proposition}
From \eqref{eq:lowerboundR} one has $R^n\gg c\lr{\xi}^{-2\nu}I$ and $\lr{\xi}^{\rho}_{\ell}R^n\gg c\lr{\xi}^{\rho-2\nu}I$ with some $c>0$ then repeating the same arguments proving \eqref{eq:rhoRn} above there is $c>0$ such that
\begin{equation}
\label{eq:Rsei}
\begin{split}
&({\rm Op}(R^n)\,v,v)\geq c\,\|\lr{D}_{\ell}^{-\nu}\,v\|^2,\\
& ({\rm Op}(\lr{\xi}^{\rho}_{\ell}R^n)\,v,v)\geq c\,\|\lr{D}_{\ell}^{-\nu+\rho/2}\,v\|^2
 \end{split}
\end{equation}
for $b\geq b_0$.
In particular ${\rm Op}(\lr{\xi}^{\rho}_{\ell}R^n)$ is nonnegative hence 
\begin{equation}
\label{eq:enkyo}
2\big|({\rm Op}(\lr{\xi}_{\ell}^{\rho}R^n)\,v,\,w)\big|\leq \delta ({\rm Op}(\lr{\xi}_{\ell}^{\rho}R^n)\,v,\,v)+\delta^{-1}({\rm Op}(\lr{\xi}_{\ell}^{\rho}R^n)\,w,\,w)
\end{equation}
for any $\delta>0$.

\subsection{Estimate of (II)}

Consider the term ${\mathsf{Re}}\,(R^nW^{n}\,U^n,\,W^{n}\,G^n\,U^n)$. 
Recall $G^n=\chi_{2h}(iA(nk,x,D)+B(nk,x))\chi_{2h}$ and with $t_n=nk$
\[
W^n\#A(t,x,\xi)\#W^{-n}=H(t,x,\xi,at_n)+R(t,at_n),\quad R(t,at_n)\in {\tilde S}^{m^*}
\]
where 
\[
W^{-n}\ :=\
\big(W_h(at_n,\xi)\big)^{-n}
\ =\
e^{-({\bar\tau}-at_n)\lr{\xi}^{\rho}_{\ell}\chi_h}\,.
\]
Then thanks to Proposition \ref{pro:b_1},
\begin{align*}
W^n\,&G^n\,W^{-n}=\chi_{2h}\,W^n(iA(t_n,x,D)+B(t_n,x))W^{-n}\,\chi_{2h}\\
&=\chi_{2h}(iH(t_n,x,D,at_n)+R(t_n,at_n))\chi_{2h}+\chi_{2h}\,W^nB(t_n,x)W^{-n}\,\chi_{2h}.
\end{align*}
Since $\chi_{2h}\in {\tilde S}^0$ and $H(t_n,x,\xi,at_n)\in {\tilde S}^1$, one sees that 
\[
\chi_{2h}\#(iH(t_n,x,\xi,at_n))\#\chi_{2h}
\ =\
i\chi_{2h}^2H(t_n,x,\xi,at_n)+{\tilde R}_n
\]
where ${\tilde R}_n\in {\tilde S}^0$ uniformly in all parameters 
satisfying
$at_n=ank\leq {\bar\tau}$. 
Define $K^n:={\tilde R}_n+\chi_{2h}\#R(t_n,at_n)\#\chi_{2h}+\chi_{2h}\#W^n\#B(t_n)\#W^{-n}\#\chi_{2h}$.   Then 
\[
W^n\#G^n\#W^{-n}\ =\
i\chi_{2h}^2H(t_n,x,\xi,at_n)\ +\
K^n(x,\xi)
\]
so,
\begin{equation}
\label{eq:omGom}
W^n\#G^n=(i\chi_{2h}^2H(t_n,x,\xi,at_n)+K^n)\#W^n\,.
\end{equation}
In addition, 
\[
K^n\in {\tilde S}^{{\bar m}},\quad
{\rm with}
\quad {\bar m}=\max\{0,m^*\}.
\]
Note that $2\nu+{\bar m}\leq \rho$ since $2\nu+{\bar m}\leq \rho$ and $\rho> 3\nu$ by \eqref{eq:rhoandnu}. Recall
\[
R_h=b\int_0^{\infty}\lr{\xi}_{\ell}^{\rho}(e^{sM^h})^*e^{sM^h}ds,\quad M^h=i\chi_{2h}^2 H(t,x,\xi,\tau)-b\,\lr{\xi}^{\rho}_{\ell}
\]
and $R^n_h=R_h(t_n,x,\xi,at_n)$ so that from \eqref{eq:RiH} it follows that
\begin{equation}
\label{eq:RH+HR}
\begin{split}
R^n(i\chi_{2h}^2H(t_n,x,\xi,at_n))+(i\chi_{2h}^2    H  &  (t_n,x,\xi,at_n))^*R^n
\\
&=-b\,\lr{\xi}_{\ell}^{\rho}+2b\,\lr{\xi}_{\ell}^{\rho}R^n.
\end{split}
\end{equation}
In view of \eqref{eq:omGom}, denoting $H(t_n)=H(t_n,x,\xi,at_n)$, one has 
\begin{align*}
2{\mathsf{Re}}\,(R^nW^{n}\, & U^n,\,W^{n}\,G^n\,U^n)=2{\mathsf{Re}}\,(W^n\,U^n,\,R^nW^n\,G^n\,U^n)
\\
&=2{\mathsf{Re}}\,(W^nU^n,\,R^n\,{\rm Op}(i\chi_{2h}^2H(t_n)+K^n)\,W^n\,U^n)
\\
&=2{\mathsf{Re}}\,(R^n\,{\rm Op}(i\chi_{2h}^2H(t_n)+K^n)W^n\,U^n,\,W^n\,U^n)
\\
&=({\rm Op}(F)W^n\,U^n,\,W^n\,U^n)\,.
\end{align*}
It follows from \eqref{eq:RH+HR} that
\begin{align*}
F&=R^n\#(i\chi_{2h}^2H(t_n)+K^n)+(i\chi_{2h}^2H(t_n)+K^n)^*\#R^n\\
&=-b\, \lr{\xi}^{\rho}_{\ell} +2b\, \lr{\xi}^{\rho}_{\ell}R^n+L^n+{\tilde L}^n\,,
\end{align*}
where $b\, L^n\in {\tilde S}^{1-\rho+3\nu}_{\rho-\nu,1-\rho+\nu}$ and ${\tilde L}^n\in {\tilde S}^{2\nu+{\bar m}}_{\rho-\nu,1-\rho+\nu}\subset {\tilde S}^{\rho}_{\rho-\nu,1-\rho+\nu}$. 
Since $\rho\geq 1-\rho+3\nu$ taking another $b_0$ if necessary one concludes
\begin{align*}
-b\,(\lr{D}_{\ell}^{\rho}\,W^n\,U^n,\,W^n\,U^n)+{\mathsf{Re}}({\rm Op}(L^n+{\tilde L}^n & )\,W^n\,U^n,\,W^n\,U^n)
\\
&\leq -\frac{b}{2}\|\lr{D}_{\ell}^{\rho/2}\,W^n\,U^n\|^2
\end{align*}
for $b\geq b_0$.  Thanks to \eqref{eq:enkyo} one has
\[
2b\, ({\rm Op}(\lr{\xi}^{\rho}_{\ell}R^n)\,W^n\,U^n,\,W^n\,U^n)\leq 4\,b \sum_{j=0}^1({\rm Op}((\lr{\xi}^{\rho}_{\ell}R^n)\,W^n\,u^{n+j},\,W^n\,u^{n+j}).
\]
Combining these estimates one obtains for $b\geq b_0$,
\begin{equation}
\label{eq:RnU}
\begin{split}
2\,{\mathsf{Re}}\,(R^nW^n\,U^n,\,W^n& \,G^n\,U^n)\leq -\frac{b}{2}\,\|\lr{D}_{\ell}^{\rho/2}\,W^n\,U^n\|^2
\\
&+4\,b \sum_{j=0}^1({\rm Op}((\lr{\xi}^{\rho}_{\ell}R^n)\,W^n\,u^{n+j},\,W^n\,u^{n+j}).
\end{split}
\end{equation}
Next study ${\mathsf{Re}}\,(R^nW^{n+1}\,G^n\,U^n,\,W^{n+1}\,U^n)$.
\begin{lemma}
\label{lem:kankan}One has
\[
W^{n+1}\#G^n\#W^{-(n+1)}=i\chi_{2h}^2H(t_n)+K^n+T^n,\quad
{\rm with}
\quad T^n\in {\tilde S}^{0}.
\]
\end{lemma}
{\bf Proof.}
Write $W^{n+1}\#G^n\#W^{-(n+1)}=W\#\big(W^n\#G^n\#W^{-n})\#W^{-1}$ so that
\[
W^{n+1}\#G^n\#W^{-(n+1)}=W\#\big(i\chi_{2h}^2H(t_n)+K^n)\#W^{-1}.
\]
Since $W^{\pm 1}\in {\tilde S}^{0}$ and $H(t_n)\in {\tilde S}^1$  it is clear that
\[
W\#(i\chi_{2h}^2H(t_n)+K^n+T^n)\#W^{-1}=i\chi_{2h}^2H(t_n)+K^n+T^n, \quad T^n\in {\tilde S}^{0}\,.
\]
This proves the lemma.
\hfill
\qed
\vskip.2cm
Lemma \ref{lem:kankan} implies that
\begin{align*}
2\,{\mathsf{Re}}\,(R^n\,W^{n+1}&G^n\,U^n,\,W^{n+1}\,U^n)\\
&=2\,{\mathsf{Re}}\,\big((R^n\,{\rm Op}(i\chi_{2h}^2H(t_n)+K^n+T^n)\,W^{n+1}U^n,\,W^{n+1}U^n)\\
&=({\rm Op}(F)W^{n+1}U^n,\,W^{n+1}U^n)
\end{align*}
with
\[
F\ :=\
R^n\#(i\chi_{2h}^2H(t_n)+K^n+T^n)+(i\chi_{2h}^2H(t_n)+K^n+T^n)^*\#R^n.
\]
Since $R^n\#T^n+(T^n)^*\#R^n\in {\tilde S}^{2\nu}_{\rho-\nu,1-\rho+\nu}$ and $\rho \geq 4\nu$ by \eqref{eq:rhonu} repeating the same arguments proving \eqref{eq:RnU} one obtains for $b\geq b_0$
\begin{equation}
\label{eq:Rn+1U}
\begin{split}
2\,{\mathsf{Re}}\,(R^nW^{n+1}\,G^n\,&U^n,\,W^{n+1}\,U^n)\leq -\frac{b}{2}\,\|\lr{D}_{\ell}^{\rho/2}\,W^{n+1}\,U^n\|^2
\\
&+4\,b \sum_{j=0}^1({\rm Op}((\lr{\xi}^{\rho}_{\ell}R^n)\,W^{n+1}\,u^{n+j},\,W^{n+1}\,u^{n+j}).
\end{split}
\end{equation}
Equations \eqref{eq:RnU} and \eqref{eq:Rn+1U} yield the following lemma.
\begin{lemma}
\label{lem:RnGU}There exist $b_0>0$ and $\ell_0>0$ such that for $b\geq b_0$ and $\ell\geq \ell_0$ one has
\begin{align*}
\frac{1}{4}\sum_{j=0}^1\,{\mathsf{Re}}\,( R^n\,W^{n+j} \,U^n,\, &W^{n+j}    \,G^n\,U^n)
\leq -\frac{b}{16}\sum_{j=0}^1\|\lr{D}_{\ell}^{\rho/2}\,W^{n+j}\,U^n\|^2
\\
&+\frac{b}{2} \sum_{i=0}^1\sum_{j=0}^1\big({\rm Op}(\lr{\xi}^{\rho}_{\ell}R^n)\,W^{n+i}\,u^{n+j},\,W^{n+i}\,u^{n+j}\big).
\end{align*}
\end{lemma}
Next  estimate $\sum_{i=0}^1{\mathsf{Re}}\,(W^{n+i}R^nW^{n+i}f^n,\,U^n)$.
Since $R^n\in {\tilde S}^{2\nu}_{\rho-\nu,1-\rho+\nu}$,
 it follows that
\begin{align*}
\big|\sum_{i=0}^1\big(W^{n+i}R^n & W^{n+i}f^n,\, U^n\big)\big|
\\
&\leq \sum_{i=0}^1\|\lr{D}_{\ell}^{-\rho/2}\,R^nW^{n+i}f^n\|\|\lr{D}_{\ell}^{\rho/2}\,W^{n+i}U^n\|\\
&\leq \frac{b}{16}\sum_{i=0}^1\|\lr{D}_{\ell}^{\rho/2}\,W^{n+i}\,U^n\|^2
+\frac{C}{b}\sum_{i=0}^1\|\lr{D}_{\ell}^{2\nu-\rho/2}\,W^{n+i}\,f^n\|^2
\,.
\end{align*}
Equation \eqref{eq:rhonu} implies that $-\nu>2\nu-\rho/2$ so
\begin{equation}
\label{eq:RnUnfn}
\begin{split}
\frac{1}{2}\sum_{i=0}^1{\mathsf{Re}}\,&( W^{n+i} R^nW^{n+i}f^n,\,U^n)
\\
&\leq \frac{b}{32}\sum_{i=0}^1\|\lr{D}_{\ell}^{\rho/2}\,W^{n+i}\,U^n\|^2
+\frac{C}{b}\sum_{i=0}^1\|\lr{D}_{\ell}^{-\nu}\,W^{n+i}\,f^n\|^2.
\end{split}
\end{equation}
Lemma \ref{lem:RnGU} together with  \eqref{eq:subCrNi} and \eqref{eq:RnUnfn}
yield the following proposition.
\begin{proposition}
\label{pro:waruiko}There exist $C>0$, $b_0>0$ and $\ell_0>0$ such that for $b\geq b_0$ and $\ell\geq \ell_0$ one has
\begin{align*}
(II)\ \leq\
 -\frac{b}{32}\sum_{i=0}^1\|\lr{D}_{\ell}^{\rho/2}\,& W^{n+i}\,U^n\|^2
 \\
+\frac{b}{2} \sum_{i=0}^1 &\sum_{j=0}^1\big({\rm Op}(\lr{\xi}^{\rho}_{\ell}R^n)\,W^{n+i}\,u^{n+j},\,W^{n+i}\,u^{n+j}\big)
\\
&+\frac{C}{b}\sum_{i=0}^1\|\lr{D}_{\ell}^{-\nu}\,W^{n+i}\,f^n\|^2.
\end{align*}
\end{proposition}
%

\subsection{Proof of Theorem \ref{thm:stability}}

First choose $b={\bar b}$ and $\ell_1$ such that Propositions \ref{pro:iiko} and \ref{pro:waruiko} and \eqref{eq:Rsei} hold with $b={\bar b}$ and $\ell\geq \ell_1$. Next choose $a={\bar a}$ such that $c\,{\bar a}\geq {\bar b}/2$ then taking \eqref{eq:Rsei} into account it follows from Propositions \ref{pro:iiko} and \ref{pro:waruiko}  that
\begin{equation}
\label{eq:I+II}
\begin{split}
(I)&+(II)
\leq  -c\,{\bar a  }\sum_{i=0}^1\sum_{j=0}^1\|\lr{D}_{\ell}^{-\nu+\rho/2}\,W^{n+i}\,u^{n+j}\|^2
\\
&-c'\,{\bar b}\sum_{i=0}^1\|\lr{D}_{\ell}^{\rho/2}\,W^{n+i}\,U^n\|^2
+C{\bar b}^{-1}\sum_{i=0}^1\|\lr{D}_{\ell}^{-\nu}\,W^{n+i}\,f^n\|^2.\end{split}
\end{equation}
%
Finally we estimate (III).  Thanks to Lemma \ref{lem:Rn+1Rn}
one has 
\begin{equation}
\label{eq:RminusR}
\begin{split}
|(III)|
\ &=\
\Big|\frac{
((R^{n+1}-R^n)W^{n+1}
u^{n+1},\,W^{n+1}u^{n+1})}{k}\Big|\\
\ &\leq\  C'\,{\bar b}^{-1}\|\lr{D}_{\ell}^{-\nu+\rho/2}\,W^{n+1}\,u^{n+1}\|^2.
\end{split}
\end{equation}

Increase ${\bar a}$  if necessary so that $c\,{\bar a}\geq 2C'\,{\bar b}^{-1}$, in view of \eqref{eq:I+II} and \eqref{eq:RminusR}, recalling \eqref{eq:I+II+III},  we conclude that
\begin{equation}
\label{eq:wa}
\begin{split}
\delta_k(R^nW^nu^n,\,W^nu^n)\leq -\frac{c}{2}{\bar a}\sum_{i=0}^1\sum_{j=0}^1\|\lr{D}_{\ell}^{-\nu+\rho/2}\,W^{n+i}\,u^{n+j}\|^2\\
-c'\,{\bar b}\sum_{i=0}^1\|\lr{D}_{\ell}^{\rho/2}\,W^{n+i}\,U^n\|^2+C\,{\bar b}^{-1}\sum_{i=0}^1\|\lr{D}_{\ell}^{-\nu}\,W^{n+i}\,f^n\|^2.
\end{split}
\end{equation}
Noting  \eqref{eq:constraint:a} and \eqref{eq:constraint:b}  we set
\[
\ell_2\ :=\ 
\max{\{{\bar a}^{6/\rho},\, {\bar b}^{1-\rho},\, \ell_1\}}\,.
\]
In what follows we assume $\ell\geq \ell_2$. Taking \eqref{eq:constraint:1} into account 
define
\[
{\bar \beta}
\ :=\
\min{\{1/2\,{\bar C},\,\log{2}/3\,{\bar a}\}}.
\]
Note that $\|\lr{D}_{\ell}^{-\nu}\,W^{n+1}\,f^n\| \leq \|\lr{D}_{\ell}^{-\nu}\,W^{n}\,f^n\|$ thanks to \eqref{eq:W_h}. Summing  \eqref{eq:wa} from $0$ to  $n-1$ yields
 \begin{equation}
 \label{eq:stab:zero}
 \begin{split}
  (R^nW^nu^n,\,W^n & u^n)+k\,\,\frac{c}{2}\,{\bar a}\sum_{p=0}^{n}\|\lr{D}_{\ell}^{-\nu+\rho/2}W^p u^p\|
  \\
 & \leq\  (R\,W^0u^0,\, W^0u^0)+C\,k\sum_{p=0}^{n-1}\|\lr{D}_{\ell}^{-\nu}\,W^{p}\,f^p\|^2.
  \end{split}
  \end{equation}
Since $W^p=e^{({\bar\tau}-{\bar a}t_p)\lr{D}_{\ell}^{\rho}\chi_h}$ 
with  $\chi_h=1$ on ${\rm supp}\,\chi_{2h}$, and recalling \eqref{eq:smallspec}, it follows from \eqref{eq:Rsei} and \eqref{eq:stab:zero} that
  \begin{align*}
  \|\lr{D}_{\ell}^{-\nu}\,e^{({\bar\tau}-{\bar a}t_n)\lr{D}_{\ell}^{\rho}}\,  u^n\|^2 & +c\,k\,{\bar a}\sum_{p=0}^n\|\lr{D}_{\ell}^{-\nu+\rho/2}e^{({\bar\tau}-{\bar a}t_p)\lr{D}_{\ell}^{\rho}}u^p\|^2
  \\
  \leq C\|\lr{D}_{\ell}^{\nu}\,& e^{{\bar\tau}\lr{D}_{\ell}^{\rho}}u^0\|^2
  +C\,k\sum_{p=0}^{ n-1}\|\lr{D}_{\ell}^{-\nu}\,e^{({\bar\tau}-{\bar a}t_p)\lr{D}_{\ell}^{\rho}}f^p\|^2.
  \end{align*}
  Equation \eqref{eq:rhonu} implies that 
$\rho/2-\nu>2\nu$  yielding the following proposition.
  \begin{proposition}
  \label{pro:CrNi_stability} There exist ${\bar \tau}>0$, ${\bar a}>0$, ${\bar \beta}>0$, $C>0$ and ${\bar \ell}\,(\geq \ell_2 )$ such that one has
\begin{equation}
\label{eq:stab:variant}
\begin{split}
\| & \lr{D}_{\ell}^{-\nu}e^{({\bar\tau}-{\bar a}t_n)\lr{D}_{\ell}^{\rho}}u^n\|^2+k\,{\bar a}\sum_{p=0}^n\|\lr{D}_{\ell}^{2\nu}e^{({\bar\tau}-{\bar a}t_p)\lr{D}_{\ell}^{\rho}}u^p\|^2
\\
&\leq C\|\lr{D}_{\ell}^{\nu}e^{{\bar\tau}\lr{D}_{\ell}^{\rho}}g\|^2
+C\,k\,\sum_{p=0}^{ n-1}\|\lr{D}_{\ell}^{-\nu}e^{({\bar\tau}-{\bar a}t_p)\lr{D}_{\ell}^{\rho}}f^p\|^2
\\
&\leq C\|\lr{D}_{\ell}^{-\nu}e^{{\bar\tau}\lr{D}_{\ell}^{\rho}}g\|^2
+C\,({\bar\tau}/{\bar a})\sup_{0\leq p\leq n-1}\|\lr{D}_{\ell}^{-\nu}e^{({\bar\tau}-{\bar a}t_p)\lr{D}_{\ell}^{\rho}}f^p\|^2
\end{split}
\end{equation}  
  for any $n\in \NN$, $k>0, \ell>0, h>0$ satisfying $nk\leq {\bar \tau}/{\bar a}$, $kh^{-1}\leq {\bar\beta}$ and $h^{-1}\geq \ell\geq {\bar \ell}$.
  \end{proposition}
   \begin{remark}
  \label{rem:jyuyo}\rm  To obtain Proposition \ref{pro:CrNi_stability} the spectral condition $\chi_{h}u^n=u^n$ is assumed while for $f^n$ no spectral condition is assumed.
  \end{remark}

\noindent
{\bf Proof of Theorem \ref{thm:stability}}: Fix $\ell={\bar\ell}$ in Proposition \ref{pro:CrNi_stability}. Since
\begin{equation}
\label{eq:xiellxi}
\lr{\xi}^{\rho}\leq \lr{\xi}^{\rho}_{{\bar\ell}}\leq {\bar\ell}^{\rho}+\lr{\xi}^{\rho},\quad \lr{\xi}\leq \lr{\xi}_{{\bar\ell}}\leq {\bar\ell}\,\lr{\xi}
\end{equation}
the proof is immediate.
\hfill
\qed

\section{Error estimates for   the spectral Crank-Nicholson scheme}
\label{sec:convergence}

\subsection{Continuous case revisited}

Start by extending  estimates \eqref{eq:conti_estimate:2} in Corollary \ref{cor:fixell} to $\partial_t^ju$ for $j=1,2$. It is clear that one can assume ${\bar\tau}\leq T$ and ${\bar a}\geq {\hat c}$. Then it is easy to examine that  Corollary \ref{cor:fixell} holds with $T={\bar\tau}$ and ${\hat c}={\bar a}$. Suppose $\partial_tu=Gu$.  Write
\[
\lr{D}_{\ell}^{\mu}G\lr{D}_{\ell}^{-\mu}=G+B_{\mu}
\]
so $\lr{D}_{\ell}^{\mu}u$ satisfies $\partial_t(\lr{D}_{\ell}^{\mu}u)=(G+B_{\mu})\lr{D}_{\ell}^{\mu}u$.  The $B_\mu$ satisfy the following bounds.
\begin{lemma}
\label{lem:Bmu}There is $A>0$ such that for any $\alpha, \beta\in\NN^d$ one has
\[
|\partial_{\xi}^{\alpha}\partial_x^{\beta}B_{\mu}(x,\xi)|\leq C_{\alpha}A^{|\beta|}|\beta|!^s\lr{\xi}_{\ell}^{-|\alpha|}\lr{x}^{-2d}.
\]
\end{lemma}
{\bf Proof}. Up to a multiplicative constant $B_{\mu}$ is given by
\begin{align*}
B_{\mu}(x,\xi)=\sum_{|\gamma|=1}\int e^{-iy\eta}\partial_{\eta}^{\gamma}\big(\lr{\xi+\eta/2}^{\mu}\lr{\xi-\eta/2}^{-\mu}\big)dyd\eta\\
\times \int \partial_x^{\gamma}G(x+\theta y,\xi)d\theta.
\end{align*}
Therefore $\partial_{\xi}^{\alpha}\partial_x^{\beta}B_{\mu}$ is, after change of variables $x+\theta y\mapsto y$,  $\theta^{-1}\eta\mapsto \eta$,  a sum of  terms
\[
\int e^{ix\eta}\partial_{\xi}^{\alpha'+\gamma}\big(\lr{\xi+\theta\eta/2}^{\mu}\lr{\xi-\theta\eta/2}^{-\mu}\big)dyd\eta
 \int e^{-iy\eta}\partial_{\xi}^{\alpha''}\partial_x^{\gamma+\beta}G( y,\xi)d\theta
\]
with $\alpha'+\alpha''=\alpha$. 
Recall that
\[
\Big|\int e^{-iy\eta}\partial_{\xi}^{\alpha''}\partial_x^{\gamma+\beta}G( y,\xi)d\theta\Big|\leq C_{\alpha''}\lr{\xi}_{\ell}^{1-|\alpha''|}A^{|\beta|}|\beta|!^se^{-c\lr{\eta}^{\rho}}
\]
with some $c>0$ (see \cite[Lemma 6.2]{CNR}). 
In addition,  it is easy to see that
\[
\big|\partial_{\eta}^{\delta}\partial_{\xi}^{\alpha'+\gamma}\big(\lr{\xi+\theta\eta/2}^{\mu}\lr{\xi-\theta\eta/2}^{-\mu}\big)\big|\leq C_{\alpha' \delta}\lr{\xi}_{\ell}^{-1-|\alpha'|}\lr{\eta}^{2|\mu|+|\delta|+1+|\alpha'|}.
\]
Using $\lr{x}^{2d}e^{ix\eta}=\lr{D_{\eta}}^{2d}e^{ix\eta}$, an
 integration by parts in $\eta$ proves the assertion.
\hfill
\qed
\vskip.2cm

Thanks to Lemma \ref{lem:Bmu} it follows from the proof of
Proposition \ref{pro:b_1}  that 
\[
e^{({\bar\tau}-{\bar a}t)\lr{\xi}_{\ell}^{\rho}}\#B_{\mu}\#e^{-({\bar\tau}-{\bar a}t)\lr{\xi}_{\ell}^{\rho}}
\ \in\ {\tilde S}^0\,.
\]
Apply Corollary \ref{cor:fixell} to  $\partial_tv=(G+B_{\mu})v$ with $v=\lr{D}_{\ell}^{\mu}u$ 
to find that 
choosing a smaller ${\bar\tau}>0$ and  larger ${\bar a}>0$ and $\ell_0$  if necessary,
\begin{equation}
\label{eq:higherD}
\|\lr{D}_{\ell}^{-\nu+\mu}e^{({\bar\tau}-{\bar a}t)\lr{D}_{\ell}^{\rho}}u(t)\|\leq C
\|\lr{D}_{\ell}^{\nu+\mu}e^{{\bar\tau}\lr{D}_{\ell}^{\rho}}u(0)\|
\end{equation}
for $0\leq t\leq {\bar\tau}/{\bar a}$ and $\ell\geq \ell_0$. Indeed in the proof of Proposition \ref{pro:CNR:1} the term $B$ satisfies $e^{(T-{\hat c}t)\lr{\xi}_{\ell}^{\rho}}\#B\# e^{(T-{\hat c}t)\lr{\xi}_{\ell}^{\rho}}\in {\tilde S}^0$
 so choosing ${\hat c}$ large, it
  is irrelevant. Write
\[
\lr{D}_{\ell}^{\mu}e^{({\bar\tau}-{\bar a}t)\lr{D}_{\ell}^{\rho}}\partial_tu
=\Big(e^{({\bar\tau}-{\bar a}t)\lr{D}_{\ell}^{\rho}}\,(G+B_{\mu})\,e^{-({\bar\tau}-{\bar a}t)\lr{D}_{\ell}^{\rho}}\Big)\lr{D}_{\ell}^{\mu}e^{({\bar\tau}-{\bar a}t)\lr{D}_{\ell}^{\rho}}u.
\]
Proposition \ref{pro:b_1} and Lemma \ref{lem:Bmu} imply  that
$e^{({\bar\tau}-{\bar a}t)\lr{\xi}_{\ell}^{\rho}}\#(G+B_{\mu})\#e^{-({\bar\tau}-{\bar a}t)\lr{\xi}_{\ell}^{\rho}}\in {\tilde S}^1$.
It  follows that
\begin{align*}
\|\lr{D}_{\ell}^{\mu}e^{({\bar\tau}-{\bar a}t)\lr{D}_{\ell}^{\rho}}\partial_tu(t)\|
\ &\leq\ C'\|\lr{D}_{\ell}^{1+\mu}e^{({\bar\tau}-{\bar a}t)\lr{D}_{\ell}^{\rho}}u(t)\|
\\
\ &\leq\  C'C\|\lr{D}_{\ell}^{1+2\nu+\mu}e^{{\bar\tau}\lr{D}_{\ell}^{\rho}}u(0)\|
\end{align*}
from  \eqref{eq:higherD}. Next assume that $A_j(t,x)$ and $B(t,x)$ are $C^1$ in time uniformly on compact sets with 
values in $G^{s'}(\RR^d)$. Since $\partial_t^2u=(\partial_tG)u+G\partial_tu$ repeating the same arguments one has
\begin{align*}
\|\lr{D}_{\ell}^{\mu} & e^{({\bar\tau}-{\bar a}t)\lr{D}_{\ell}^{\rho}}\partial_t^2u(t)\| \\
\ &\leq\ 
C''\big(\|\lr{D}_{\ell}^{1+\mu}e^{({\bar\tau}-{\bar a}t)\lr{D}_{\ell}^{\rho}}u(t)\|
\ +\ \|\lr{D}_{\ell}^{1+\mu}e^{({\bar\tau}-{\bar a}t)\lr{D}_{\ell}^{\rho}}\partial_tu(t)\|\big)
\\
\ &\leq\
 C'''\|\lr{D}_{\ell}^{2+2\nu+\mu}e^{{\bar\tau}\lr{D}_{\ell}^{\rho}}u(0)\|.
\end{align*}
Choosing $\mu=-\nu+i$, $i=0,1,2$ one obtains the following lemma.
\begin{lemma}
\label{lem:C2est} Assume that $A_j(t,x)$  and $B(t,x)$ are $C^1$
 in time uniformly on compact sets with 
values in $G^{s'}(\RR^d)$ and that $\partial_tu=Gu$. Then there exist $C>0$, $\ell_0>0$ such that
\[
\|\lr{D}_{\ell}^{-\nu+i}e^{({\bar\tau}-{\bar a}t)\lr{D}_{\ell}^{\rho}}\partial_t^ju(t)\|
\ \leq\
 C\|\lr{D}_{\ell}^{i+j+\nu}e^{{\bar\tau}\lr{D}_{\ell}^{\rho}}u(0)\|
\]
for $0\leq t\leq {\bar \tau}/{\bar a}$, $\ell\geq \ell_0$ and $0\le i,j\le 2$.
\end{lemma}

 \subsection{Error  estimate for the spectral Crank-Nicholson scheme}
  

Suppose that $u(t,x)$
  satisfies
\begin{equation}
\label{eq:Exsol}
\partial_tu(t,x)=G(t,x,D)u(t,x)
\end{equation}
where $G(t,x,D)=iA(t,x,D)+B(t,x)$. 
 Denote ${\tilde u}=\chi_{2h} u$  so that $\chi_{h}{\tilde u}={\tilde u}$. Thus 
 \begin{equation}
\label{eq:Exsol_bis}
\partial_t{\tilde u}=G{\tilde u}+f,\quad f=[\chi_{2h},G]u.
\end{equation}
 %
 
 Next estimate to what extent ${\tilde u}(t_n,x)$ satisfies the difference scheme. The error, denoted by $g(n)=g(n,\cdot)$, is given by
\[
\frac{{\tilde u}(t_{n+1})-{\tilde u}(t_n)}{k}
\ -\
G^n\,\frac{{\tilde u}(t_{n+1})+{\tilde u}(t_n)}{2}
\ :=\ g(n)
\]
where $G^n=\chi_{2h}(iA(nk,x,D)+B(nk,x))\chi_{2h}$. 
Note that
\begin{equation}
\label{eq:suppg}
{\rm supp}\,{\mathcal F}(g(n))
\ \subset\
 {\rm supp}\,\chi_{2h}(\cdot).
\end{equation}
The approximate solution $u^n=u^n_h$ satisfies
\[
\frac{u^{n+1}-u^n}{k}-G^n\,\frac{u^{n+1}+u^n}{2}=0.
\]
At $t=0$ the approximate solution is equal to the spectral
truncation of the exact solution, $u^0=\chi_{2h}g={\tilde u}(0)$.  

Noting  ${\rm supp}{\mathcal F}\big({\tilde u}(t_n)-u^n\big)\subset {\rm supp}\,\chi_{2h}$ and hence $\chi_h({\tilde u}(t_n)-u^n)={\tilde u}(t_n)-u^n$,
 Proposition \ref{pro:CrNi_stability} implies
\begin{equation}
\label{eq:gosa}
\|\lr{D}_{\ell}^{-\nu}\,W^n\,({\tilde u}(t_{n})-u^n)\|^2
\ \leq\
 C\,k
\sum_{l=0}^{ n-1}
\|\lr{D}_{\ell}^{-\nu}\,W^l\,g(l)\|^2
\end{equation}
for any $t_n=kn\leq {\bar \tau}/{\bar a}$.
\begin{lemma}
\label{lem:Estg}
There is $C>0$ so that
\begin{align*}
\|\lr{D}_{\ell}^{-\nu}W^jg(j)\|
\leq
 C\,(k+h)\|\lr{D}^{2+\nu}_{\ell}e^{{\bar\tau}\lr{D}^{\rho}_{\ell}}u(0)\|
 \end{align*}
for $0\leq j\leq n-1$ and $0\leq t_n\leq {\bar\tau}/{\bar a}$.
\end{lemma}
{\bf Proof.} Use \eqref{eq:Exsol_bis} to write
\[
g(j) =g(j)-\big({\tilde u}_t(t_j)
-
G(t_j){\tilde u}(t_j)-f(j)\big).
\]
The triangle inequality yields
\begin{equation}
\label{eq:Wjgj}
\begin{split}
\|\lr{D}_{\ell}^{-\nu}&W^jg(j)\|\ \leq\ \Big\|\lr{D}_{\ell}^{-\nu}W^j\Big(\frac{{\tilde u}(t_{j+1})-{\tilde u}(t_j)}{k}  - {\tilde u}_t(t_j)\Big)\Big\|
\\
&+
\Big\|\lr{D}_{\ell}^{-\nu}W^j\Big(G^j\Big(\frac{{\tilde u}(t_{j+1})+{\tilde u}(t_j)}{2}-{\tilde u}(t_j)\Big)\Big)\Big\|
\\
&+\|\lr{D}_{\ell}^{-\nu}W^j\big(G(t_j)-G^j\big)\,{\tilde u}(t_j)\|
+\|\lr{D}_{\ell}^{-\nu}W^jf(j)\|.
\end{split}
\end{equation}
Write
\[
\lr{D}_{\ell}^{-\nu}W^j\Big(\frac{{\tilde u}(t_{j+1})-{\tilde u}(t_j)}{k}  - {\tilde u}_t(t_j)\Big)=\frac{1}{k}\int_{t_j}^{t_{j+1}}ds\int_{t_j}^s\lr{D}_{\ell}^{-\nu}W^j\partial_t^2{\tilde u}(s')ds'
\]
and note that
\[
W^j\partial_t^2{\tilde u}(s')=e^{{\bar a}(s'-t_j)\lr{D}_{\ell}^{\rho}\chi_h}e^{({\bar\tau}-{\bar a}s')\lr{D}_{\ell}^{\rho}\chi_h}\partial_t^2{\tilde u}(s').
\]
Since $0\leq s'-t_j\leq k$ if $t_j\leq s'\leq t_{j+1}$ it follows from \eqref{eq:W_h} that
\begin{align*}
\|\lr{D}_{\ell}^{-\nu}&W^j\partial_t^2{\tilde u}(s')\|\leq 2\|\lr{D}_{\ell}^{-\nu}e^{({\bar\tau}-{\bar a}s')\lr{D}_{\ell}^{\rho}\chi_h}\partial_t^2{\tilde u}(s')\|
\\
&\leq 2\|\lr{D}_{\ell}^{-\nu}e^{({\bar\tau}-{\bar a}s')\lr{D}_{\ell}^{\rho}}\partial_t^2{ u}(s')\|\leq C\|\lr{D}_{\ell}^{2+\nu}e^{{\bar\tau}\lr{D}_{\ell}^{\rho}}u(0)\|
\end{align*}
thanks to Lemma \ref{lem:C2est}. Therefore one has
\[
\Big\|\lr{D}_{\ell}^{-\nu}W^j\Big(\frac{{\tilde u}(t_{j+1})-{\tilde u}(t_j)}{k}  - {\tilde u}_t(t_j)\Big)\Big\|\leq C\,k\,\|\lr{D}_{\ell}^{2+\nu}e^{{\bar\tau}\lr{D}_{\ell}^{\rho}}u(0)\|.
\]
Turn to the second term on the right-hand side of \eqref{eq:Wjgj}. 
Use
\[
\lr{D}_{\ell}^{-\nu}W^j\Big(G^j\Big(\frac{{\tilde u}(t_{j+1})+{\tilde u}(t_j)}{2}-{\tilde u}(t_j)\Big)\Big)=\frac{1}{2}\int_{t_j}^{t_{j+1}}\lr{D}_{\ell}^{-\nu}W^jG^j\partial_t{\tilde u}(s')ds'
\]
to write
\[
\lr{D}_{\ell}^{-\nu}W^jG^j=\lr{D}_{\ell}^{-\nu}W^jG^jW^{-j}\big(W^{j}e^{-({\bar\tau}-{\bar a}s')\lr{D}_{\ell}^{\rho}\chi_h}\big)e^{({\bar\tau}-{\bar a}s')\lr{D}_{\ell}^{\rho}\chi_h}.
\]

Proposition \ref{pro:b_1}  implies that 
 $\lr{\xi}_{\ell}^{-\nu}\#W^j\#G^j\#W^{-j}\in {\tilde S}^{1-\nu}$.
 In addition, $W^je^{-({\bar\tau}-{\bar a}s')\lr{D}_{\ell}^{\rho}\chi_h}=e^{{\bar a}(s'-t_j)\lr{D}_{\ell}^{\rho}\chi_h}$ when $0\leq s'-t_j\leq k$.
Repeat the same arguments as above to find
\[
\|\lr{D}_{\ell}^{-\nu}W^jG^j\partial_t{\tilde u}(s')\|
\ \leq\
 C\|\lr{D}_{\ell}^{2+\nu}e^{{\bar\tau}\lr{D}_{\ell}^{\rho}}u(0)\|
\quad
{\rm 
for}
\quad t_j\leq s'\leq t_{j+1}\,.
\]
 Then 
\[
\Big\|\lr{D}_{\ell}^{-\nu}W^j\Big(G^j\Big(\frac{{\tilde u}(t_{j+1})+{\tilde u}(t_j)}{2}-{\tilde u}(t_j)\Big)\Big)\Big\|\leq C\,k\,\|\lr{D}_{\ell}^{2+\nu}e^{{\bar\tau}\lr{D}_{\ell}^{\rho}}u(0)\|.
\]

Next study the third and fourth term on the right-hand side of \eqref{eq:Wjgj}.
\begin{lemma}
\label{lem:uuM}
Let $\alpha\geq 0$. There is $C>0$ such that
\[
\|(I-\chi_{2 h})u\|\leq C\,h^{\alpha}\|\lr{D}_{\ell}^{\alpha}u\|.
\]
\end{lemma}
{\bf Proof.}
Since $1-\chi_{2h}(\xi)=0$ unless $|\xi|\geq h^{-1}$ one has
\begin{align*}
\|(I-\chi_{2 h})u\|^2
&=\int(1-\chi_{2 h}(\xi))^2\lr{\xi}^{-2\alpha}_{\ell}\lr{\xi}_{\ell}^{2\alpha}|{\hat u}(\xi)|^2d\xi
\\
&\leq Ch^{2\alpha}\int \lr{\xi}^{2\alpha}_{\ell}|{\hat u}(\xi)|^2d\xi
= C\Big(h^{\alpha}\|\lr{D}_{\ell}^{\alpha}u\|\Big)^2
\end{align*}
which proves the assertion.
\hfill
\qed
\vskip.2cm

Since $G^j-G(t_j)=\chi_{2h}G(t_j)(\chi_{2h}-I)+(\chi_{2h}-I)G(t_j)$ one can write
\begin{align*}
\lr{D}_{\ell}^{-\nu}W^j(G^j-G(t_j))=\chi_{2h}\big(\lr{D}_{\ell}^{-\nu}W^jG(t_j)W^{-j}\big)(\chi_{2h}-I)W^j\\
+(\chi_{2h}-I)\big(\lr{D}_{\ell}^{-\nu}W^jG(t_j)W^{-j}\big)W^j.
\end{align*}
Using $\lr{\xi}_{\ell}^{-\nu}\#W^j\#G(t_j)\#W^{-j}\in {\tilde S}^{1-\nu}$ 
together with Lemma \ref{lem:uuM} one finds
\begin{align*}
\| & \lr{D}_{\ell}^{-\nu}W^j(G^j-G(t_j)){\tilde u}(t_j)\|
\\
&\le C\,\|\lr{D}_{\ell}^{1-\nu}(\chi_{2h}-I)W^j{\tilde u}(t_j)\|
+C\,h\,\|\lr{D}_{\ell}^{1-\nu}W^jG(t_j)W^{-j}W^j{\tilde u}(t_j)
\|
\\
&\leq C'\,h\,\|\lr{D}_{\ell}^{2-\nu}W^j{\tilde u}(t_j)\|
\leq C'\,h\,\|\lr{D}_{\ell}^{2-\nu}e^{({\bar\tau}-{\bar a}t_j)\lr{D}_{\ell}^{\rho}}u(t_j)\|.
\end{align*}
Therefore by Lemma \ref{lem:C2est},
\[
\|\lr{D}_{\ell}^{-\nu}W^j(G^j-G(t_j)){\tilde u}(t_j)\|\leq C\,h\,\|\lr{D}_{\ell}^{2+\nu}e^{{\bar\tau}\lr{D}_{\ell}^{\rho}}u(0)\|.
\]
Turn to $f(j):=[\chi_{2h},G(t_j)]u(t_j)$. Since 
\[
[\chi_{2h},G(t_j)]=\chi_{2h}G(t_j)(I-\chi_{2h})-(I-\chi_{2h})G(t_j)\chi_{2h}
\]
 repeating the same arguments as above one obtains that
\[
\|\lr{D}_{\ell}^{-\nu}W^jf(j)\|\leq C\,h\,\|\lr{D}_{\ell}^{2+\nu}e^{{\bar\tau}\lr{D}_{\ell}^{\rho}}u(0)\|.
\]
This finishes the proof of Lemma \ref{lem:Estg}. 
\hfill
\qed
\vskip.2cm
\subsection{Proof of Theorem \ref{thm:convergence}}

Noting that ${\rm supp}{\mathcal F}\big({\tilde u}(t_n)-u^n\big)\subset {\rm supp}\,\chi_{2h}$ and $\chi_h=1$ on the support of $\chi_{2h}$ 
 it follows from \eqref{eq:gosa} and Lemma \ref{lem:Estg} that
\begin{equation}
\label{eq:untilu}
\begin{split}
\|\lr{D}_{\ell}^{-\nu}\,e^{({\bar\tau}-{\bar a}t_n)\lr{D}_{\ell}^{\rho}}\,&({\tilde u}(t_{n})-u^n)\|
\\
&\leq \  C\,\sqrt{{\bar\tau}/{\bar a}}\,(k+h)\|\lr{D}^{2+\nu}_{\ell}e^{{\bar\tau}\lr{D}^{\rho}_{\ell}}u(0)\|.
\end{split}
\end{equation}
Since $\lr{\xi}_{\ell}\leq \sqrt{3}h^{-1}$ on the support of $\chi_{2h}$,
\eqref{eq:untilu} implies that 
\begin{equation}
\label{eq:untilu:bis}
\|e^{({\bar\tau}-{\bar a}t_n)\lr{D}_{\ell}^{\rho}}\,({\tilde u}(t_{n})-u^n)\|
\leq  C\,\sqrt{{\bar\tau}/{\bar a}}\,(k+h)h^{-\nu}\|\lr{D}^{2+\nu}_{\ell}e^{{\bar\tau}\lr{D}^{\rho}_{\ell}}u(0)\|.
\end{equation}
Finally estimate $\|\lr{D}_{\ell}^{-\nu}W^n(u(t_n)-{\tilde u}(t_n))\|$. 
Since $u(t_n)-{\tilde u}(t_n)=(1-\chi_{2h})u(t_n)$ the same arguments as above prove that
\begin{equation}
\label{eq:utilu}
\|\lr{D}_{\ell}^{-\nu}\,e^{({\bar\tau}-{\bar a}t_n)\lr{D}_{\ell}^{\rho}}\,(u(t_n)-{\tilde u}(t_n))\|
\ \leq\  Ch^2\|\lr{D}_{\ell}^{2+\nu}e^{{\bar\tau}\lr{D}_{\ell}^{\rho}}u(0)\|.
\end{equation}
Similarly one has
\begin{equation}
\label{eq:utilu:bis}
\|e^{({\bar\tau}-{\bar a}t_n)\lr{D}_{\ell}^{\rho}}\,(u(t_n)-{\tilde u}(t_n))\|
\ \leq\ 
 Ch^{2-\nu}\|\lr{D}_{\ell}^{2+\nu}e^{{\bar\tau}\lr{D}_{\ell}^{\rho}}u(0)\|.
\end{equation}
Combining \eqref{eq:untilu}, \eqref{eq:untilu:bis} and \eqref{eq:utilu}, \eqref{eq:utilu:bis} yields
the following proposition.
\begin{proposition}
\label{pro:CrNi_accuracy}There exist ${\bar\tau}>0, {\bar a}>0, {\bar\beta }>0, C>0$ and ${\bar \ell}>0$ such that for any  exact solution $u$ to \eqref{eq:Exsol} with Cauchy data $u(0)$ such that $\lr{D}^{2+\nu}_{\ell}e^{{\bar \tau}\lr{D}_{\ell}^{\rho}}u(0)\in L^2$ one has 
\[
\|\lr{D}_{\ell}^{-\nu}e^{({\bar\tau}-{\bar a}t_n)\lr{D}_{\ell
}^{\rho}}(u(t_n)-u^n)\|
\leq C\,(k+h)\|\lr{D}_{\ell}^{2+\nu}e^{{\bar\tau}\lr{D}_{\ell}^{\rho}}u(0)\|
\]
and
\[
\|e^{({\bar\tau}-{\bar a}t_n)\lr{D}_{\ell
}^{\rho}}(u(t_n)-u^n)\|
\leq C\,(k+h)h^{-\nu}\|\lr{D}_{\ell}^{2+\nu}e^{{\bar\tau}\lr{D}_{\ell}^{\rho}}u(0)\|
\]
for any $0\leq t_n=nk\leq {\bar\tau}/{\bar a}$, $kh^{-1}\leq {\bar\beta}$ and $h^{-1}\geq \ell\geq {\bar \ell}$.
\end{proposition}
%

\begin{remark}\rm
      \label{rem:cfl} 
      In order for  a difference approximation to be accurate,
      the time discretization must be taken sufficiently fine
      \cite{CFL}.  Here Proposition \ref{pro:CrNi_accuracy} shows that one could constrain $k$ to satisfy a CFL type condition  $kh^{-1}\leq {\bar\beta}$. More precisely,
       the proof shows that it suffices to constrain $k$ to satisfy
\[
kh^{-1}\leq 1/2{\bar C},\quad kh^{-\rho}\leq \log{2}/3{\bar a}.
\]
  \end{remark}

\noindent
{\bf Proof of Theorem \ref{thm:convergence}}: Taking \eqref{eq:xiellxi} into account it is enough to choose $\ell={\bar\ell}$ in Proposition \ref{pro:CrNi_accuracy}.
\hfill
\qed



\begin{thebibliography}{99}

\bibitem{Bronsh:1} 
M.D.Bronshtein, { The Cauchy problem for hyperbolic operators with characteristics of variable multiplicity}, Trudy Moskov. Mat. Obsc. {\bf 41} (1980), 83-99; English translation: Trans. Moscow Math. Soc. {\bf 41} (1982), 87-103.






\bibitem{CNR}
F.Colombini, T.Nishitani, J.Rauch,
{ Weakly hyperbolic systems by symmetrization},
Ann. Scuola Norm. Sup. Pisa, {\bf 19} (2019),  217-251.





\bibitem{CR2} 
F.Colombini, J.Rauch,  { Numerical analysis of very weakly well-posed hyperbolic Cauchy problems}, IMA J. Numer. Anal. {\bf 35} (2015), 989-1010.

\bibitem{CFL}
R.Courant, K.O.Friedrichs, H.Lewy,
{ \"Uber die partiellen differenzengleichungen der mathematischen physik}, Math. Ann.
{\bf 100} (1928), 32-74.

\bibitem{garding}
L. G\aa rding, 
Linear hyperbolic partial differential equations with constant coefficients. 
Acta Math. {\bf 85}(1951), 1-62. 



\bibitem{Hobook} 
L.H\"ormander, {\it  The Analysis of Linear Partial Differential Operators. III. Pseudo-Differential Operators}, Springer, Berlin, 1994.




\bibitem{Lax2}
P.D.Lax, 
Asymptotic  solutions of oscillatory initial value problems,
Duke Math. J. {\bf 24}(1957), 627-646.


\bibitem{Lax1}
P.D.Lax, { Differential equations, difference equations and matrix theory}, Comm. Pure Appl. Math., {\bf 11} (1958), 175-194.

\bibitem{mizohata}
S. Mizohata, 
 Some remarks on the Cauchy problem. J. Math. Kyoto Univ. {\bf 1} (1961/62), 109-127.

\bibitem{nishitani}
T. Nishitani, On the Lax-Mizohata theorem in the analytic and Gevrey classes. J. Math. Kyoto Univ. {\bf 18} (1978),  509-521.

\bibitem{sabrina}  S.Petit-Bergez,  { Probl\`emes faiblement bien pos\'es:
discr\'etisations et applications},
Th\`ese Docteur de Math\'ematiques Universit\'e de Paris 13, 2006.


  

\bibitem{YaNo}
M.Yamaguti, T.Nogi, {An algebra of pseudo difference schemes and its applications}, Publ. Res. Inst. Math. Sci. {\bf 3} (1967), 151-166.














\end{thebibliography}
\end{document}